\documentclass[11pt,reqno]{amsart}

\usepackage{lmodern}
\usepackage{amsmath}
\usepackage{amssymb}
\usepackage{amscd}
\usepackage{mathrsfs}
\usepackage{mathtools}
\usepackage{booktabs}
\usepackage{microtype}
\usepackage[hidelinks]{hyperref}
\hypersetup{
  pdftitle={Finite diagonalizable torsors and Kummer cohomology over semiring schemes},
  pdfauthor={Chandrasekhar Gokavarapu, Sajani Lavanya Madasi, Rajeev Muthu, Sekhar Babu Gosala, Naveen Kumar Kakumanu},
  pdfkeywords={semiring scheme; diagonalizable torsor; strong grading; derived cohomology; Kummer gerbe; fpqc descent}
}

\newtheorem{theorem}{Theorem}[section]
\newtheorem{lemma}[theorem]{Lemma}
\newtheorem{proposition}[theorem]{Proposition}
\newtheorem{corollary}[theorem]{Corollary}

\theoremstyle{definition}
\newtheorem{definition}[theorem]{Definition}

\theoremstyle{remark}
\newtheorem{remark}[theorem]{Remark}

\DeclareMathOperator{\Pic}{Pic}
\DeclareMathOperator{\Spec}{Spec}

\DeclareMathOperator{\Hom}{Hom}
\DeclareMathOperator{\Isom}{Isom}
\DeclareMathOperator{\Tors}{Tors}
\DeclareMathOperator{\Ext}{Ext}
\DeclareMathOperator{\Br}{Br}
\DeclareMathOperator{\coker}{coker}

\newcommand{\NN}{\mathbb N}

\newcommand{\TT}{\mathbb T}
\newcommand{\Gm}{\mathbf G_{\mathrm m}}
\newcommand{\OO}{\mathcal O}
\newcommand{\Kum}{\mathcal K}

\newcommand{\fpqc}{\mathrm{fpqc}}

\newcommand{\ot}{\mathbin{\otimes}}

\newcommand{\longiso}{\xrightarrow{\;\sim\;}}

\title[Finite diagonalizable torsors]{Finite diagonalizable torsors and Kummer cohomology over
semiring schemes}

\author{Chandrasekhar Gokavarapu}
\address{Department of Mathematics\\ Acharya Nagarjuna University,
Guntur, Andhra Pradesh 522510, India}
\curraddr{Department of Mathematics\\ Government College (Autonomous),
Rajahmundry, Andhra Pradesh 533105, India}
\email{chandrasekhargokavarapu@gcrjy.ac.in}
\email{chandrasekhargokavarapu@gmail.com}
\urladdr{https://orcid.org/0009-0006-5306-371X}
\thanks{Corresponding author: Chandrasekhar Gokavarapu.}

\author{Sajani Lavanya Madasi}
\address{Department of Mathematics\\ Government College (Autonomous),
Rajahmundry, Andhra Pradesh 533105, India}
\email{sajanimaths@gcrjy.ac.in}

\author{Rajeev Muthu}
\address{Department of Mathematics\\ Government College (Autonomous),
Rajahmundry, Andhra Pradesh 533105, India}
\email{euclid.gdc@gmail.com}
\urladdr{https://orcid.org/0009-0000-5549-1108}

\author{Sekhar Babu Gosala}
\address{Department of Mathematics\\ Government College (Autonomous),
Rajahmundry, Andhra Pradesh 533105, India}
\email{sekhar1983@gcrjy.ac.in}
\urladdr{https://orcid.org/0009-0009-5410-9267}

\author{Naveen Kumar Kakumanu}
\address{Department of Mathematics\\ Government College (Autonomous),
Rajahmundry, Andhra Pradesh 533105, India}
\email{ramanawinmaths@gmail.com}
\urladdr{https://orcid.org/0000-0001-5329-3343}

\subjclass[2020]{Primary 14F20, 14L15; Secondary 14A20, 16W50, 16Y60}
\keywords{semiring scheme, diagonalizable torsor, strongly graded
semialgebra, Kummer cohomology, fpqc descent, root gerbe}
\date{}

\begin{document}
\raggedbottom

\begin{abstract}
Classical descent treats a finite diagonalizable torsor as a
higher-rank locally free object.  Over semirings that route is
incomplete: fpqc-local freeness is not known to imply Zariski-local
freeness in arbitrary rank.  We show that diagonalizable symmetry
bypasses this obstruction.  The coordinate semialgebra splits into
character pieces, each piece descends in rank one, and the torsor
identity makes their multiplication maps invertible.  This gives a
natural equivalence between fpqc \(D_X(Q)\)-torsors and Picard-strong
\(Q\)-gradings, and represents every torsor by a finite locally free
morphism of rank \(|Q|\).

The same viewpoint proves that finite-index monomial maps of split
tori are finite locally free torsors and yields matrix Kummer
classification sequences.  A concrete example is the nontrivial
Boolean torsor
\(\Spec\mathbb B[v^{\pm1}]\to\Spec\mathbb B[u^{\pm1}]\),
\(u\mapsto v^m\), whose unit-power classes contribute
\(\mathbb Z/m\mathbb Z\) to degree-one cohomology.  A two-term
resolution of the character group then gives a
presentation-independent universal-coefficient filtration for
diagonalizable fpqc cohomology.  This formalism organizes, but does not
compute, \(\mathbf G_{\mathrm m}\)-cohomology.  In degree two it
identifies the Kummer boundary of a line bundle with its root gerbe.
On projective space over the Boolean or real tropical semifield all
\(\mu_m\)-torsors vanish, whereas the root gerbe of \(\mathcal O(1)\)
has exact order \(m\).  Thus torsors and gerbes retain information that
can disappear under ordinary ring completion.
\end{abstract}

\maketitle
\setcounter{tocdepth}{2}
\tableofcontents

\section{Introduction}
\label{sec:introduction}

\subsection{A motivating puzzle}

Consider the Boolean semifield \(\mathbb B=\{0,1\}\), with
\(1+1=1\), and fix \(m\geq2\).  The group of \(\mathbb B\)-valued
points of \(\mu_m\) is trivial: the only invertible element of
\(\mathbb B\) is \(1\).  It is therefore tempting to expect that
\(\mu_m\)-torsors over simple Boolean schemes should also be trivial.
The affine Boolean torus gives the opposite answer.  The map
\[
 \operatorname{Spec}\mathbb B[v^{\pm1}]
 \longrightarrow
 \operatorname{Spec}\mathbb B[u^{\pm1}],
 \qquad u\longmapsto v^m,
\]
is a finite free \(\mu_m\)-torsor, but it has no section: a section
would produce a Laurent monomial whose \(m\)th power is \(u\), which
is impossible.  Thus a group scheme can have no visible nonidentity
base-valued points and still possess geometrically nontrivial torsors.

This elementary example contains the three themes of the paper.  The
coordinate semialgebra splits into the \(m\) character pieces
\[
 \mathbb B[u^{\pm1}]\oplus
 \mathbb B[u^{\pm1}]v\oplus\cdots\oplus
 \mathbb B[u^{\pm1}]v^{m-1};
\]
adjoining roots is a finite torus cover; and the obstruction to a
section is a Kummer class.  Our purpose is to show that the same
picture works for every finite split diagonalizable group over a
semiring scheme, and then to follow it one degree higher from torsors
to root gerbes.

\subsection{Why ordinary descent is not enough}

Let \(Q\) be a finite abelian group and let \(D(Q)\) denote the
diagonalizable group with character group \(Q\).  Over an ordinary
scheme, a \(D(Q)\)-action is encoded by a \(Q\)-grading, and Kummer
theory extracts from a presentation of \(Q\) a finite flat torus
isogeny.  These two descriptions connect diagonalizable torsors,
strongly graded algebras, abelian covers, and root constructions
\cite{GrothendieckSGA3,DemazureGabriel1970,Waterhouse1979,
ChaseHarrisonRosenberg1965,Dade1980,NastasescuVanOystaeyen2004,
Pardini1991,Cadman2007,Milne1980,StacksKummer,StacksRoot}.

Transferring this picture to semiring schemes creates a descent
obstruction.  Semimodule flatness is categorical, quotients are
controlled by congruences, and the usual ring-theoretic equivalences
among finite projectivity, flatness, and local freeness are not
available in arbitrary rank
\cite{Golan1999,HebischWeinert1998,BorgerJun2025}.  In particular,
Borger and Jun prove the needed equivalence in rank one, while the
corresponding higher-rank assertion remains open in their framework.
Consequently, merely descending the regular representation of
\(D(Q)\), which has rank \(|Q|\), does not prove that the resulting
semimodule is Zariski-locally free.

\subsection{The character-by-character escape}

The main result resolves this obstruction for finite split
diagonalizable torsors.  The coaction splits the coordinate
semialgebra into its character components.  Each component descends
separately from a rank-one free semimodule, so the rank-one theorem
applies; the torsor identity then forces every homogeneous
multiplication map to be invertible.  Thus this particular
higher-rank object is reconstructed from finitely many effective
rank-one descent problems.  We make no claim about arbitrary
higher-rank semimodules.

Informally, one should think of the characters of \(Q\) as colors.
The coordinate algebra of the trivial torsor has one free rank-one
piece of each color.  Descent preserves the colors, so instead of
descending one mysterious object of rank \(|Q|\), we descend \(|Q|\)
line-like objects separately.  The torsor identity says that
multiplying a piece of color \(p\) by one of color \(q\) gives the
piece of color \(p+q\).  This is precisely the extra structure that
turns the rank-one pieces back into a torsor.

The remaining results develop consequences of that geometric input.
Once a finite-index monomial map of split tori is proved to be a
finite locally free fpqc epimorphism over semirings, a free resolution
of \(Q\) becomes a two-term resolution of \(D(Q)\) by tori.  This
expresses the derived fpqc cohomology of \(D(Q)\) functorially in
terms of the generally noncomputed cohomology of \(\Gm\).  The resulting
universal-coefficient filtration shows that degree one and degree
two measure different failures of divisibility in the unit and
Picard groups.  In the projective examples below, the resulting
torsors vanish while specific root gerbes survive.

\subsection{Three tempting shortcuts and why they fail}

The route to these conclusions is easier to see by separating three
plausible arguments that do not suffice.

\begin{enumerate}
\item \emph{Descend the regular representation all at once.}
An fpqc trivialization does descend the underlying rank-\(|Q|\)
semimodule.  What it does not presently supply over an arbitrary
semiring is a theorem turning that fpqc form into a Zariski-locally
free semimodule.  The character decomposition is therefore not merely
notation: it replaces this unavailable higher-rank step by finitely
many established rank-one steps.

\item \emph{Test surjectivity on base-valued points.}
Over \(\mathbb B\), the only \(m\)th root of unity is \(1\), and a
typical unit need not have an \(m\)th root in the base.  Nevertheless,
the power morphism is fpqc-surjective because roots appear after the
finite free extension constructed in Section~\ref{sec:power}.
Pointwise surjectivity over the original semiring is thus the wrong
test for a torsor statement.

\item \emph{Pass first to the universal ring.}
For a nonzero idempotent semiring, additive group completion is the
zero ring.  Ring completion consequently erases both the affine
Boolean torsor and the projective Boolean root gerbe studied below.
It can provide a useful shadow when nondegenerate, but it cannot be a
conservative proof method in characteristic one.
\end{enumerate}

These failures also explain the vertical organization of the paper.
The rank-one argument establishes the geometric object, the finite
torus cover supplies the exact character resolution, and only then is
ordinary derived sheaf cohomology applied.  None of the three layers
can be substituted for the preceding one.

\subsection{Main theorem: graded reconstruction}

Write \(\operatorname{PicGr}_Q(X)\) for the groupoid of systems
\[
 \bigl((\mathcal L_q)_{q\in Q},u,(m_{p,q})_{p,q\in Q}\bigr),
\]
where every \(\mathcal L_q\) is a line bundle,
\(u\colon\mathcal O_X\!\xrightarrow{\sim}\mathcal L_0\), and
\[
 m_{p,q}\colon
 \mathcal L_p\otimes\mathcal L_q
 \xrightarrow{\sim}\mathcal L_{p+q}
\]
satisfy the unit, associativity, and commutativity diagrams.

\begin{theorem}[Graded reconstruction]
\label{thm:intro-graded}
For every semiring scheme \(X\) and every finite abelian group \(Q\),
there is a natural equivalence
\[
 \Tors_{\fpqc}(X,D_X(Q))
 \xrightarrow{\ \sim\ }
 \operatorname{PicGr}_Q(X).
\]
Under this equivalence a torsor has coordinate semialgebra
\[
 \mathcal B=\bigoplus_{q\in Q}\mathcal L_q.
\]
It is represented by a finite locally free morphism of rank
\(|Q|\).  Conversely, every Picard-strong \(Q\)-grading reconstructs
a \(D_X(Q)\)-torsor through the Hopf--Galois isomorphism.
\end{theorem}

This combines Theorems~\ref{thm:graded-torsor-equivalence}
and~\ref{thm:global-graded-descent}.  The conclusion is stronger than
fpqc-local representability: it identifies the internal summands of
the descended coordinate algebra and proves their multiplication
maps are isomorphisms.  The assignment
\[
 q\longmapsto[\mathcal L_q]
\]
is therefore a homomorphism \(Q\to\Pic(X)\), but the grading retains
the coherent multiplication data which this homomorphism forgets.
The proof is specific to diagonalizable coactions; it does not claim
higher-rank local freeness for arbitrary semimodules.

Although the corresponding statement over rings is classical, the
argument is not obtained by citing the usual strong-grading theorem
verbatim.  The classical theorem starts with an already represented
graded algebra and works in an additive module category.  Here one
starts with an fpqc sheaf torsor.  One must first descend its affine
coordinate semialgebra, prove that the coaction decomposition survives
effective descent coefficientwise, turn each fpqc rank-one summand into
a Zariski line bundle, and only then descend the homogeneous
multiplication isomorphisms.  The last two steps use, respectively, the
rank-one theorem of Borger--Jun and fpqc conservativity.  This isolates
the precise extra input required beyond Dade's ring-theoretic
criterion.

\subsection{Kummer cohomology as a consequence}

The geometric theorem has a cohomological continuation.  Here
\(R\Gamma_{\fpqc}\) is derived global sections in the abelian category
of sheaves of abelian groups on the fpqc site.  It is unrelated to a
hypothetical derived category of semimodules.

\begin{theorem}[Relative derived character theorem]
\label{thm:intro-derived}
For every morphism \(f\colon X\to Y\) of semiring schemes and finite
abelian group \(Q\), there is a natural isomorphism in the derived
category of abelian fpqc sheaves on \(Y\):
\[
 Rf_*D_X(Q)
 \xrightarrow{\ \sim\ }
 R\mathcal{H}\!om_{\mathbb Z,Y}
 \bigl(Q_Y,Rf_*\Gm{}_{X}\bigr).
\]
Its cohomology sheaves carry a natural two-step
universal-coefficient filtration.  In particular, taking derived
global sections, for every \(i\geq1\) there is a
presentation-independent exact sequence
\[
 0\longrightarrow
 \Ext^1_{\mathbb Z}
 \bigl(Q,H^{i-1}_{\fpqc}(X,\Gm)\bigr)
 \longrightarrow
 H^i_{\fpqc}(X,D_X(Q))
\]
\[
 \longrightarrow
 \Hom_{\mathbb Z}
 \bigl(Q,H^i_{\fpqc}(X,\Gm)\bigr)
 \longrightarrow0.
\]
The sequence need not split naturally.
\end{theorem}

The relative statement and its global edge are proved in
Theorems~\ref{thm:relative-derived-character-formula},
\ref{thm:relative-universal-coefficient}, and
\ref{thm:universal-coefficient}.  They follow from a length-one free
resolution with equal-rank terms,
\[
 0\longrightarrow\mathbb Z^r\xrightarrow{\,M\,}\mathbb Z^r
 \longrightarrow Q\longrightarrow0
\]
and the fpqc-exact sequence
\[
 1\longrightarrow D_X(Q)\longrightarrow\Gm^r
 \xrightarrow{\varphi_M}\Gm^r\longrightarrow1.
\]
Thus the matrix Kummer calculation is the degree-one edge of a
single derived statement.  The formula does not by itself evaluate
the \(\Gm\)-cohomology of \(X\); its role is to reduce
diagonalizable cohomology to that input in a presentation-independent
way.

In degree one, using
\(H^0(X,\Gm)=\Gamma(X,\mathcal O_X^\times)\) and
\(H^1(X,\Gm)=\Pic(X)\), one obtains
\[
 0\longrightarrow
 \Ext^1_{\mathbb Z}
 \bigl(Q,\Gamma(X,\mathcal O_X^\times)\bigr)
 \longrightarrow H^1_{\fpqc}(X,D_X(Q))
\]
\[
 \longrightarrow
 \Hom_{\mathbb Z}(Q,\Pic(X))
 \longrightarrow0.
\]
For \(Q=\coker M\), this is the matrix Kummer sequence.
The presentation remains useful geometrically: the corresponding
torus isogeny is finite free of rank \(|\det M|\), its homotopy fiber
on the Picard groupoid classifies the torsors, and Smith normal form
recovers the invariant cyclic factors.

The unit term is already nonzero on a characteristic-one affine
scheme.  For
\(X=\operatorname{Spec}\mathbb B[u^{\pm1}]\), the morphism
\[
 \operatorname{Spec}\mathbb B[v^{\pm1}]
 \longrightarrow X,
 \qquad u\longmapsto v^m,
\]
is a nontrivial finite free \(\mu_m\)-torsor of rank \(m\).  Since the
units of \(\mathbb B[u^{\pm1}]\) are precisely the monomials \(u^a\),
the classes of the root semialgebras \(z^m=u^a\) inject
\(\mathbb Z/m\mathbb Z\) into \(H^1_{\fpqc}(X,\mu_m)\).  This gives an
explicit affine application that does not depend on a prior Picard
group computation; see Section~\ref{sec:affine-boolean-torsor}.

\subsection{The degree-two edge and projective examples}

Put
\[
 \Br'_{\fpqc}(X)=H^2_{\fpqc}(X,\Gm).
\]
We make no comparison here with Azumaya semialgebras.  Abelian
second cohomology classifies banded gerbes
\cite{Giraud1971,StacksSecondCohomologyGerbes}, and the degree-two
edge of the derived character theorem gives the following.

\begin{theorem}[Diagonalizable gerbes and root detection]
\label{thm:intro-gerbes}
There is a natural exact sequence
\[
\begin{aligned}
 0&\longrightarrow \Ext^1_{\mathbb Z}(Q,\Pic(X))
 \longrightarrow H^2_{\fpqc}(X,D_X(Q))\\
 &\longrightarrow \Hom_{\mathbb Z}(Q,\Br'_{\fpqc}(X))
 \longrightarrow0.
\end{aligned}
\]
For \(Q=\mathbb Z/m\mathbb Z\), its left-hand map is
\[
 \Pic(X)/m\Pic(X)\lhook\joinrel\longrightarrow
 H^2_{\fpqc}(X,\mu_m);
\]
the class of \(L\) is the gerbe
\(\sqrt[m]{L/X}\) of \(m\)th roots of \(L\).  This gerbe is neutral
exactly when \([L]\in m\Pic(X)\), and its cohomological order equals
the order of \([L]\) modulo \(m\Pic(X)\).
\end{theorem}

This result is Theorems~\ref{thm:diagonalizable-gerbe-sequence}
and~\ref{thm:root-gerbe-detection}.  It exposes a vertical
distinction that degree-one theory alone suppresses: torsion in
\(\Pic(X)\) contributes to torsors, whereas failure of divisibility
of \(\Pic(X)\) contributes to gerbes.

For a totally ordered idempotent semifield \(S\), Jun's computation
\(\Pic(\mathbb P^d_S)\cong\mathbb Z\)
\cite{Jun2017} makes the distinction explicit.  If
\(S=\mathbb T_{\mathbb R}\) or \(S=\mathbb B\), then for every
\(m\geq2\)
\[
 H^1_{\fpqc}(\mathbb P^d_S,\mu_m)=0,
 \qquad
 [\sqrt[m]{\mathcal O(1)/\mathbb P^d_S}]
 \in H^2_{\fpqc}(\mathbb P^d_S,\mu_m)
\]
has exact order \(m\).  Over the Boolean semifield, additive group
completion is the zero ring, so this nonzero class has no
nondegenerate classical shadow.  It provides a concrete test of the
difference between the degree-one and degree-two edges.

\subsection{Relation to previous work, scope, and organization}

Over ordinary rings, diagonalizable actions and strong gradings, the
Hopf--Galois criterion, finite abelian covers, and Kummer sequences are
classical; see
\cite{GrothendieckSGA3,DemazureGabriel1970,Waterhouse1979,
ChaseHarrisonRosenberg1965,Dade1980,NastasescuVanOystaeyen2004,
Pardini1991,Milne1980}.  We do not claim these principles as new.  The
new point of Theorem~\ref{thm:intro-graded} is their semiring-scheme
form together with representability and finite local freeness: the
character summands permit use of the rank-one descent theorem of
Borger and Jun \cite{BorgerJun2025}, whereas descent of the full
rank-\(|Q|\) regular representation alone would require an unavailable
higher-rank implication.  Similarly, the homological algebra behind a
universal-coefficient sequence is standard.  Its use here depends on
the semiring-specific proof that the relevant monomial torus map is an
fpqc epimorphism; after that geometric step, the derived formula
organizes the torsor and gerbe consequences rather than constituting a
new general universal-coefficient principle or a method for computing
\(\Gm\)-cohomology.

Our semiring schemes are the Zariski-glued relative schemes fixed in
Section~\ref{sec:foundations}.  Their foundations lie in relative
algebraic geometry and its semiring and blueprint realizations
\cite{ToenVaquie2009,Marty2012,Durov2007,Deitmar2005,Deitmar2008,
ConnesConsani2010,Lorscheid2012,Lorscheid2017,Lorscheid2026}.
The tropical calculation uses only the projective Picard theorem in
this model.  Other approaches to tropical schemes and tropical
linear geometry provide context but are not silently identified with
it
\cite{GiansiracusaGiansiracusa2016,MaclaganSturmfels2015,
MaclaganRincon2018,JunMinchevaTolliver2019,
JunMinchevaTolliver2024,Lorscheid2022,Lorscheid2023,
FloresWeibel2014}.

The terminology ``Kummer'' also occurs in logarithmic flat geometry
and in root-stack constructions
\cite{Kato1989,Olsson2003,IllusieNakayamaTsuji2013,
AbramovichOlssonVistoli2008,BorneVistoli2012,TalpoVistoli2018,
GillibertZhao2024}.  Our topology is the ordinary fpqc topology on
semiring schemes, and our root object in degree two is a banded
gerbe.  We do not claim a logarithmic comparison, a theory for
nonsplit groups of multiplicative type, or a calculation of
\(\Br'_{\fpqc}(X)\) itself.  Standard descent and stack language is
used in the sense of
\cite{GrothendieckSGA1,Vistoli2005,LaumonMoretBailly2000,Olsson2016};
derived-functor conventions follow
\cite{Weibel1994,StacksCohomologySites}.

\subsection{A reader's map}

The paper may be read in three layers.  For the geometric mechanism,
Section~\ref{sec:foundations} fixes the descent conventions,
Section~\ref{sec:power} constructs finite-index torus covers, and
Section~\ref{sec:graded-torsors} proves graded reconstruction.  For
explicit classifications, Section~\ref{sec:matrix-kummer} gives the
matrix fiber and Section~\ref{sec:affine-boolean-torsor} returns to the
Boolean example above with complete coordinate formulas.  For the
cohomological layer, Section~\ref{sec:derived} proves the derived
character formula, Section~\ref{sec:base-change} studies classical
shadows, and Section~\ref{sec:tropical} separates torsors from gerbes
on tropical projective space.  The final section isolates four
concrete problems suggested by the method.

\section{Semiring geometry and fpqc descent}
\label{sec:foundations}

\subsection{Semimodules and flatness}

All semirings are commutative and unital.  The zero semiring is
allowed when it occurs as a pushout or a ring completion, but every
geometric base considered in the main statements is nonzero.  For a
semiring \(A\), the category \(\operatorname{Mod}(A)\) consists of
commutative additive monoids with a unital \(A\)-action.  Tensor
products are those of the closed symmetric monoidal category
\(\operatorname{Mod}(A)\); in particular,
\[
 \Hom_A(M\ot_A N,P)
 \cong
 \Hom_A\bigl(M,\Hom_A(N,P)\bigr).
\]
This is the tensor product used in relative algebraic geometry and in
the semiring-scheme literature
\cite{Golan1999,HebischWeinert1998,ToenVaquie2009,
Lorscheid2017,Jun2017,BorgerJun2025}.

We use the categorical notion of flatness.  An \(A\)-semimodule \(M\)
is \emph{flat} if
\[
 M\ot_A -\colon\operatorname{Mod}(A)\longrightarrow
 \operatorname{Mod}(A)
\]
preserves finite limits.  It is \emph{faithful} if a morphism becomes
an isomorphism after tensoring with \(M\) only when it was already an
isomorphism.  An \(A\)-algebra \(B\) is faithfully flat when its
underlying \(A\)-semimodule is flat and faithful.  This is the
definition adopted in~\cite[Section~4]{BorgerJun2025}.

\begin{lemma}
\label{lem:finite-free-faithful}
If \(r\geq1\), the free \(A\)-semimodule \(A^r\) is flat and faithful.
Consequently, an \(A\)-algebra that is free of positive finite rank as
an \(A\)-semimodule is faithfully flat.
\end{lemma}

\begin{proof}
For every \(A\)-semimodule \(M\), finite distributivity of the tensor
product gives
\[
  A^r\ot_A M\cong M^r.
\]
Finite products and finite coproducts agree in
\(\operatorname{Mod}(A)\), and limits in a finite power are computed
coordinatewise.  Hence \(M\mapsto M^r\) preserves finite limits.  If
\(f^r\colon M^r\to N^r\) is an isomorphism, then \(f\) is an
isomorphism: an inverse is obtained by inserting into one coordinate,
applying the inverse of \(f^r\), and projecting from the same
coordinate.  Thus the functor is faithful in the required sense.
\end{proof}

\begin{remark}
\label{rem:flatness-warning}
Over rings, finite projective modules are flat and are locally free in
the usual topology.  Over semirings, some classical equivalences fail
in higher rank.  The rank-one situation is substantially better:
Borger and Jun prove that invertible, Zariski-locally free of rank
one, and fpqc-locally free of rank one are equivalent
\cite[Theorem~1.1]{BorgerJun2025}.  The arguments below use finite
free modules directly for the Kummer covers and use this rank-one
theorem for line bundles.
\end{remark}

\subsection{The geometric category and the fpqc topology}

Write \(\operatorname{CAlg}(A)\) for commutative \(A\)-algebras and
\(\operatorname{Aff}_A=\operatorname{CAlg}(A)^{\mathrm{op}}\).
We use the following fixed convention throughout the paper.  An
\(A\)-semiring scheme is a sheaf on the relative Zariski site of
\(\operatorname{Aff}_A\) which admits a Zariski cover by representable
affine functors.  Equivalently, for the Zariski topology used here, it
may be described by gluing locally semiringed affine spaces
\cite[Remarks~4.2, 5.2 and~5.3]{BorgerJun2025}.
This is the relative-scheme framework of To\"en--Vaqui\'e specialized
to semimodules; see
\cite{ToenVaquie2009,Marty2012,Lorscheid2017}.  We do not transfer
statements to other notions of semiring scheme without a comparison
functor.

On affines the fpqc topology is generated by finite families
\[
 \{\Spec B_i\longrightarrow\Spec B\}_{i\in I}
\]
for which \(\prod_i B_i\) is faithfully flat over \(B\).
Faithfully flat descent for semimodules is comonadic
\cite[Theorem~4.1]{BorgerJun2025}; the corresponding assertion for
commutative algebras follows from the same descent formalism
\cite[Proposition~4.9]{BorgerJun2025}.  We shall use the following
consequences.

\begin{proposition}[Fpqc descent package]
\label{prop:descent-package}
Let \(A\to B\) be faithfully flat.
\begin{enumerate}
\item \(A\)-semimodules are equivalent to \(B\)-semimodules equipped
with effective descent data.
\item Commutative \(A\)-algebras satisfy effective descent along
\(A\to B\).
\item Finite presentation, dualizability, and invertibility of
semimodules are fpqc-local properties.
\item Invertible semimodules form an fpqc stack.
\item A morphism of semimodules is an isomorphism if and only if it
becomes an isomorphism after base change to \(B\).
\end{enumerate}
\end{proposition}

\begin{proof}
The first two assertions are
\cite[Theorem~4.1 and Proposition~4.9]{BorgerJun2025}.  Fpqc locality
of finite presentation and dualizability is proved in
\cite[Propositions~6.2 and 8.4]{BorgerJun2025}; invertibility is
\cite[Corollary~10.3]{BorgerJun2025}.  The last assertion is the
faithfulness, or conservativity, in the definition of a faithfully
flat morphism.
\end{proof}

For a quasi-coherent semimodule \(\mathcal F\) on a semiring scheme
\(X\), ``finite locally free of rank \(r\)'' always means
Zariski-locally isomorphic to \(\OO_X^r\).  We shall not use
fpqc-local finite freeness as a synonym.  This distinction is
essential: for \(r>1\), it is not known over arbitrary semirings
whether fpqc-local finite freeness implies Zariski-local finite
freeness \cite[Theorem~1.2 and the discussion following it]
{BorgerJun2025}.  Section~\ref{sec:graded-torsors} proves the stronger
Zariski-local conclusion for finite diagonalizable torsors by
decomposing their coordinate algebras into line bundles.

A \emph{line bundle} on a semiring scheme \(X\) is an invertible
\(\OO_X\)-semimodule.  We denote its dual by
\[
 L^\vee=\mathcal Hom_{\OO_X}(L,\OO_X).
\]
Evaluation gives \(L^\vee\ot L\cong\OO_X\).  The group of isomorphism
classes under tensor product is \(\Pic(X)\).  The preceding descent
theorem and the comparison of rank-one definitions ensure that this
agrees with the Zariski definition used in
\cite{Jun2017,JunMinchevaTolliver2019}.

\subsection{Monoid semirings and diagonalizable group schemes}

Let \(M\) be a commutative monoid, written additively.  Its monoid
semiring over \(A\) is
\[
 A[M]=
 \left\{
   \sum_{x\in M}^{\mathrm{finite}} a_x[X^x] : a_x\in A
 \right\},
\]
with convolution multiplication
\(X^xX^y=X^{x+y}\).  Formal sums have unique coefficients; no
cancellation or subtraction is involved.

\begin{lemma}
\label{lem:monoid-base-change}
For every morphism \(A\to A'\) of commutative semirings there is a
natural isomorphism
\[
  A'\ot_A A[M]\longiso A'[M].
\]
If \(M\) is a finite set of cardinality \(r\), then \(A[M]\) is free
of rank \(r\) as an \(A\)-semimodule.
\end{lemma}

\begin{proof}
Both sides represent the functor that assigns to an \(A'\)-algebra
\(C\) the set of monoid homomorphisms \(M\to(C,\cdot)\).  The second
assertion follows from the formal basis
\(\{X^x:x\in M\}\).
\end{proof}

If \(G\) is an abelian group, \(A[G]\) is a Hopf \(A\)-semialgebra:
\[
 \Delta(X^g)=X^g\ot X^g,\qquad
 \varepsilon(X^g)=1,\qquad
 S(X^g)=X^{-g}.
\]
It represents the diagonalizable group functor
\[
 D_A(G)(C)=\Hom_{\mathrm{Ab}}(G,C^\times)
\]
on \(A\)-algebras.  In particular,
\[
 \Gm=D_{\NN}(\mathbb Z)=\Spec\NN[\mathbb Z].
\]
This is the direct semiring analogue of the character-group
construction for diagonalizable group schemes over rings
\cite[Expos{\'e}~VIII]{GrothendieckSGA3};
see also \cite{DemazureGabriel1970,Waterhouse1979}.
For \(m\geq1\), put \(C_m=\mathbb Z/m\mathbb Z\) and define
\begin{equation}
\label{eq:mu-coordinate}
 \mu_m=D_{\NN}(C_m)=\Spec\NN[C_m].
\end{equation}
Thus, for every semiring \(A\),
\[
 \mu_m(A)=\{\zeta\in A^\times:\zeta^m=1\}.
\]
By Lemmas~\ref{lem:finite-free-faithful}
and~\ref{lem:monoid-base-change}, \(\mu_m\) is finite locally free
of rank \(m\) and faithfully flat over \(\Spec\NN\).  Notice that
this conclusion concerns the representing group scheme, not merely
its set of points.  For example, \(\mu_m(\TT_{\mathbb Z})\) is
trivial, although \(\mu_m\) has nontrivial torsors over
\(\TT_{\mathbb Z}\)-schemes; see
Section~\ref{sec:tropical}.

\subsection{Torsors and first cohomology}

Let \(G\) be an affine group scheme over a semiring scheme \(X\).  A
\(G\)-torsor on the fpqc site of \(X\) is an fpqc sheaf \(P\) with a
\(G\)-action which is fpqc-locally equivariantly isomorphic to \(G\).
We use the standard sheaf-theoretic convention for torsors and
contracted products, modeled on classical descent
\cite{GrothendieckSGA1,Giraud1971,Vistoli2005}.
If \(X=\Spec A\), \(G\) is affine, and a torsor is trivialized by an
affine fpqc cover, effective descent for commutative algebras in
Proposition~\ref{prop:descent-package} shows that the torsor is
represented by an affine \(X\)-scheme.  For a general finite locally
free \(G\), this argument yields an fpqc-locally finite free
coordinate semimodule; it does not by itself prove Zariski-local
freeness in higher rank.  For \(G=D_X(Q)\) with \(Q\) finite, the
missing stronger conclusion will follow from
Theorem~\ref{thm:graded-torsor-equivalence}.  We write
\[
 H^1_{\fpqc}(X,G)
\]
for the pointed set of isomorphism classes of torsors.  It is an
abelian group when \(G\) is commutative.  This notation refers to
torsor cohomology; it does not require an additive derived-functor
theory for sheaves of semimodules.

The standard argument with frames identifies \(\Gm\)-torsors with
line bundles.  Indeed, a line bundle has an fpqc sheaf of frames, and
an fpqc \(\Gm\)-torsor produces a line bundle by extension of
structure group.  Hence
\begin{equation}
\label{eq:H1-Gm-Pic}
 H^1_{\fpqc}(X,\Gm)\cong\Pic(X).
\end{equation}
This is compatible with the Zariski computation
\(\Pic(X)\cong H^1(X,\OO_X^\times)\) in
\cite{Jun2017,JunMinchevaTolliver2019}, because line bundles admit
Zariski local frames and all rank-one notions agree.

\section{The finite-flat power morphism}
\label{sec:power}

\subsection{Finite-index character lattices}

We first prove a diagonalizable-group statement that contains the
power morphism.  Its classical ring analogue belongs to the theory of
groups of multiplicative type
\cite[Expos{\'e}~VIII]{GrothendieckSGA3}; the point here is that the
coordinate proof is coefficientwise and remains valid over every
commutative semiring.

\begin{theorem}[Finite-index diagonalizable torsor]
\label{thm:finite-index-torsor}
Let \(A\) be a commutative semiring and let
\[
 \iota\colon\Lambda'\hookrightarrow\Lambda
\]
be an injective homomorphism of abelian groups.  Suppose that
\(Q=\Lambda/\Lambda'\) is finite of order \(d\).  The induced
morphism
\begin{equation}
\label{eq:finite-index-map}
 D_A(\Lambda)\longrightarrow D_A(\Lambda')
\end{equation}
is finite locally free of rank \(d\), faithfully flat, and a torsor
under \(D_A(Q)\) for the fpqc topology.  These assertions commute
with arbitrary base change \(A\to A'\).
\end{theorem}

\begin{proof}
Choose a set \(R\subset\Lambda\) of representatives for
\(\Lambda/\Lambda'\).  Every \(\lambda\in\Lambda\) has a unique
expression \(\lambda=\lambda'+\rho\), with
\(\lambda'\in\Lambda'\) and \(\rho\in R\).  Grouping formal
coefficients by cosets gives an isomorphism of
\(A[\Lambda']\)-semimodules
\begin{equation}
\label{eq:finite-index-free-decomposition}
 \bigoplus_{\rho\in R}A[\Lambda']e_\rho
 \longiso A[\Lambda],
 \qquad
 a_\rho e_\rho\longmapsto a_\rho X^\rho .
\end{equation}
It is free of rank \(d\), so \eqref{eq:finite-index-map} is finite
locally free and faithfully flat by
Lemma~\ref{lem:finite-free-faithful}.

The quotient sequence
\[
 0\longrightarrow\Lambda'\longrightarrow\Lambda
 \longrightarrow Q\longrightarrow0
\]
induces \(D_A(Q)\to D_A(\Lambda)\).  For an \(A\)-algebra \(C\), a
point of
\[
 D_A(\Lambda)\times_{D_A(\Lambda')}D_A(\Lambda)
\]
is a pair of characters
\(\chi_1,\chi_2\colon\Lambda\to C^\times\) with equal restrictions
to \(\Lambda'\).  Their ratio \(\chi_1\chi_2^{-1}\) is trivial on
\(\Lambda'\) and hence factors uniquely through a character
\(\bar\chi\colon Q\to C^\times\).  Therefore
\[
 D_A(Q)\times_A D_A(\Lambda)
 \longrightarrow
 D_A(\Lambda)\times_{D_A(\Lambda')}D_A(\Lambda),
 \qquad
 (\bar\chi,\chi)\longmapsto(\bar\chi\chi,\chi)
\]
is an isomorphism on \(C\)-points, naturally in \(C\).  Yoneda's
lemma gives the torsor identity.  Finally,
Lemma~\ref{lem:monoid-base-change} identifies every scalar extension
of the coordinate map with the same construction over \(A'\).
\end{proof}

\begin{corollary}[Integral-matrix torus covers]
\label{cor:matrix-isogeny}
Let \(M\colon\mathbb Z^r\to\mathbb Z^r\) be an integral matrix with
\(\det M\neq0\).  The corresponding morphism of split semiring tori
\[
 \varphi_M\colon(\Gm{}_{A})^r\longrightarrow(\Gm{}_{A})^r
\]
is a finite locally free fpqc torsor of rank \(|\det M|\), under the
finite locally free diagonalizable group scheme
\[
 D_A(\coker M).
\]
\end{corollary}

\begin{proof}
The map \(M\) is injective and its cokernel has cardinality
\(|\det M|\).  Apply Theorem~\ref{thm:finite-index-torsor}.
\end{proof}

\subsection{The one-dimensional power map}

Fix an integer \(m\geq1\).  We distinguish the Laurent generators on
the target and source of the power morphism by writing
\[
 R_A=A[s,s^{-1}]=A[\mathbb Z],
 \qquad
 S_A=A[t,t^{-1}]=A[\mathbb Z].
\]
The \(A\)-algebra morphism
\begin{equation}
\label{eq:power-coordinate-map}
 [m]^*\colon R_A\longrightarrow S_A,\qquad s\longmapsto t^m,
\end{equation}
represents the map \([m]\colon\Gm{}_{A}\to\Gm{}_{A}\) on \(A\)-schemes.

\begin{proposition}[Coset decomposition]
\label{prop:coset-decomposition}
As an \(R_A\)-semimodule through~\eqref{eq:power-coordinate-map},
\(S_A\) is free of rank \(m\), with basis
\[
 1,t,\ldots,t^{m-1}.
\]
More precisely, multiplication induces an isomorphism
\begin{equation}
\label{eq:coset-decomposition}
 \bigoplus_{r=0}^{m-1}R_A\,e_r
 \longiso S_A,
 \qquad
 a(s)e_r\longmapsto a(t^m)t^r.
\end{equation}
\end{proposition}

\begin{proof}
Every integer \(k\) has a unique expression \(k=mq+r\), with
\(q\in\mathbb Z\) and \(0\leq r<m\).  Hence a Laurent polynomial
\[
 f(t)=\sum_{k\in\mathbb Z}^{\mathrm{finite}} a_k t^k
\]
has a unique decomposition
\[
 f(t)=
 \sum_{r=0}^{m-1}
 \left(\sum_{q\in\mathbb Z}^{\mathrm{finite}}a_{mq+r}s^q\right)t^r.
\]
The uniqueness is coefficientwise in the formal monoid semiring and
does not use additive cancellation.  This is exactly the inverse of
\eqref{eq:coset-decomposition}.
\end{proof}

\begin{corollary}
\label{cor:power-finite-flat}
For every commutative semiring \(A\), the power morphism
\[
 [m]\colon\Gm{}_{A}\longrightarrow\Gm{}_{A}
\]
is finite locally free of constant rank \(m\), finitely presented,
and faithfully flat.
\end{corollary}

\begin{proof}
Proposition~\ref{prop:coset-decomposition} makes its coordinate
semialgebra free of rank \(m\).  Faithful flatness follows from
Lemma~\ref{lem:finite-free-faithful}.  The presentation by one
invertible generator with the relation \(s=t^m\) is finite.
\end{proof}

\begin{remark}
\label{rem:stronger-than-root}
The statement is stronger than saying that each unit admits an
\(m\)th root after some fpqc extension.  It gives one universal
finite locally free cover of \(\Gm\), of a fixed rank independent of
the base semiring and of the chosen unit.
\end{remark}

\subsection{Kernel and torsor identity}

The group scheme \(\mu_m\) acts on \(\Gm\) by multiplication.  At
the level of \(C\)-points, for every \(A\)-algebra \(C\), the action is
\[
 \mu_m(C)\times C^\times\longrightarrow C^\times,
 \qquad
 (\zeta,x)\longmapsto\zeta x.
\]

\begin{proposition}
\label{prop:kernel}
The group scheme \(\mu_m=\Spec\NN[C_m]\) is the scheme-theoretic
kernel of \([m]\colon\Gm\to\Gm\).
\end{proposition}

\begin{proof}
For every semiring \(C\), the kernel of
\(C^\times\to C^\times\), \(x\mapsto x^m\), is
\(\{\zeta\in C^\times:\zeta^m=1\}\).  By
\eqref{eq:mu-coordinate}, this functor is represented by \(\mu_m\).
Yoneda's lemma gives the assertion.
\end{proof}

\begin{theorem}[Universal Kummer torsor]
\label{thm:power-torsor}
The power morphism
\[
 [m]\colon\Gm\longrightarrow\Gm
\]
is a \(\mu_m\)-torsor for the fpqc topology.  It is finite locally
free of rank \(m\).  Equivalently, the canonical morphism
\begin{equation}
\label{eq:torsor-isomorphism}
 \mu_m\times\Gm
 \longrightarrow
 \Gm\times_{[m],\,\Gm,\,[m]}\Gm,
 \qquad
 (\zeta,y)\longmapsto(\zeta y,y),
\end{equation}
is an isomorphism, and \([m]\) is an fpqc cover.
\end{theorem}

\begin{proof}
The covering assertion is
Corollary~\ref{cor:power-finite-flat}.  To prove
\eqref{eq:torsor-isomorphism}, evaluate on an arbitrary semiring
\(C\).  A point of the fiber product is a pair
\((x,y)\in C^\times\times C^\times\) satisfying \(x^m=y^m\).
Then
\[
 \zeta=xy^{-1}
\]
satisfies \(\zeta^m=1\), and \((x,y)=(\zeta y,y)\).  This
construction is inverse to~\eqref{eq:torsor-isomorphism} and is
natural in \(C\).  Yoneda's lemma completes the proof.
\end{proof}

\begin{corollary}[Fpqc Kummer sequence]
\label{cor:kummer-sheaf-sequence}
The sequence
\begin{equation}
\label{eq:kummer-sheaf-sequence}
 1\longrightarrow\mu_m\longrightarrow\Gm
 \xrightarrow{[m]}\Gm\longrightarrow1
\end{equation}
is exact as a sequence of abelian fpqc sheaves on semiring schemes.
\end{corollary}

\begin{proof}
Proposition~\ref{prop:kernel} gives exactness at the first two terms.
Theorem~\ref{thm:power-torsor} makes \([m]\) an fpqc cover and hence
an epimorphism of fpqc sheaves.
\end{proof}

\subsection{Fibers and root semialgebras}

Let \(u\in A^\times\).  It defines an \(A\)-point
\[
 u\colon\Spec A\longrightarrow\Gm{}_{A}.
\]
The pullback of Theorem~\ref{thm:power-torsor} along \(u\) will be
denoted
\begin{equation}
\label{eq:root-torsor}
 \mathcal R_m(u)
 =
 \Spec A\times_{u,\,\Gm{}_{A},\,[m]}\Gm{}_{A}.
\end{equation}

\begin{proposition}[Canonical root semialgebra]
\label{prop:root-algebra}
The torsor \(\mathcal R_m(u)\) is affine, with coordinate semialgebra
\begin{equation}
\label{eq:root-algebra}
 A\langle u^{1/m}\rangle
 =
 A\ot_{A[s,s^{-1}]}A[t,t^{-1}],
\end{equation}
where \(s\mapsto u\) on the left and \(s\mapsto t^m\) on the right.
It is free of rank \(m\) over \(A\), with basis
\[
 1,z,\ldots,z^{m-1},
\]
and \(z^m=u\).  In particular \(z\) is a unit, with
\[
 z^{-1}=u^{-1}z^{m-1}.
\]
\end{proposition}

\begin{proof}
The affine fiber-product formula gives~\eqref{eq:root-algebra}.
Base-changing the free decomposition
\eqref{eq:coset-decomposition} along \(A[s,s^{-1}]\to A\) gives the
displayed \(A\)-basis.  The image \(z\) of \(t\) satisfies \(z^m=u\);
the formula for its inverse follows immediately.
\end{proof}

For later calculations it is useful to give a presentation that does
not hide the semiring congruence.  On the free \(A\)-semimodule
\[
 B_u=\bigoplus_{r=0}^{m-1}Az^r
\]
define multiplication on basis elements by
\begin{equation}
\label{eq:wrapped-multiplication}
 z^iz^j
 =
 u^{\lfloor(i+j)/m\rfloor}
 z^{\,i+j-m\lfloor(i+j)/m\rfloor}.
\end{equation}
Division of \(i+j+k\) by \(m\) proves associativity.  The resulting
semialgebra has the universal property of the congruence quotient
\[
 A[t]/(t^m\sim u).
\]
Since \(u\) is a unit, \(t\) becomes a unit, so this quotient agrees
with~\eqref{eq:root-algebra}.  This construction proves normal-form
uniqueness directly and avoids treating the relation as an ordinary
ideal quotient.

\begin{corollary}
\label{cor:root-torsor}
The morphism
\[
 \mathcal R_m(u)\longrightarrow\Spec A
\]
is a finite locally free \(\mu_{m,A}\)-torsor of rank \(m\).
It represents the functor
\[
 C\longmapsto\{z\in C^\times:z^m=u_C\}
\]
on \(A\)-algebras.
\end{corollary}

\begin{proof}
It is the base change of the torsor in
Theorem~\ref{thm:power-torsor}.  The functor-of-points description
follows from the defining fiber product.
\end{proof}

\begin{proposition}[Universality]
\label{prop:root-universal}
Let \(A\to C\) be a semiring morphism and let \(w\in C^\times\) with
\(w^m=u_C\).  There is a unique \(A\)-algebra morphism
\[
 A\langle u^{1/m}\rangle\longrightarrow C
\]
sending \(z\) to \(w\).  The construction of
\(A\langle u^{1/m}\rangle\) commutes with arbitrary base change
\(A\to A'\).
\end{proposition}

\begin{proof}
The first assertion is the universal property of the fiber product,
or equivalently of the congruence presentation above.  For the
second, use associativity of pushouts:
\[
 A'\ot_A
 \bigl(A\ot_{A[s,s^{-1}]}A[t,t^{-1}]\bigr)
 \cong
 A'\ot_{A'[s,s^{-1}]}A'[t,t^{-1}].
\]
\end{proof}

\begin{remark}
The construction remains valid when \(m\) is not invertible in the
base semiring.  No derivative criterion is used, and no claim of
etaleness is made.  Finite local freeness and the fpqc topology are
the appropriate uniform statements.
\end{remark}

\section{Strongly graded descent for diagonalizable torsors}
\label{sec:graded-torsors}

The higher-rank descent issue isolated in
Section~\ref{sec:foundations} can be resolved for diagonalizable
torsors because their coordinate algebras carry more structure than
arbitrary fpqc forms of a free semimodule.  A diagonalizable action is
a grading, and the homogeneous pieces of a torsor descend separately
as line bundles.  Over rings this principle belongs to the theory of
diagonalizable actions, strongly graded rings, and Hopf--Galois
extensions
\cite[Expos{\'e}~VIII]{GrothendieckSGA3}
\cite{Dade1980,NastasescuVanOystaeyen2004,Schauenburg2004}.
Those results explain the expected algebraic shape but do not supply
the descent statement needed here: before strong grading can be used,
an fpqc sheaf torsor must be represented and its descended coordinate
semimodule must be shown Zariski-locally free.  We prove the
coefficientwise semiring form and isolate the rank-one input that
closes this gap.

\subsection{Coactions and finite gradings}

Let \(Q\) be a finite abelian group, written additively, and let
\(A\) be a commutative semiring.  A \(Q\)-grading on an
\(A\)-semimodule \(B\) is a decomposition
\[
 B=\bigoplus_{q\in Q}B_q.
\]
It is a grading of an \(A\)-semialgebra if
\[
 1\in B_0,\qquad B_pB_q\subseteq B_{p+q}.
\]
The grading is called \emph{Picard-strong} if the structural map
\(A\to B_0\) and every multiplication map
\begin{equation}
\label{eq:strong-multiplication}
 B_p\ot_A B_q\longrightarrow B_{p+q}
\end{equation}
are isomorphisms.  In that case \(B_q\) is invertible, with inverse
\(B_{-q}\), and
\[
 B=\bigoplus_{q\in Q}B_q
\]
is Zariski-locally free of rank \(|Q|\).

\begin{lemma}[Coaction--grading correspondence]
\label{lem:coaction-grading}
Let \(B\) be an \(A\)-semimodule.  Coactions
\[
 \rho\colon B\longrightarrow B\ot_A A[Q]
\]
of the Hopf semialgebra \(A[Q]\) are naturally equivalent to
\(Q\)-gradings of \(B\).  Under this equivalence,
\[
 B_q=\{b\in B:\rho(b)=b\ot X^q\}.
\]
Coactions on commutative \(A\)-semialgebras correspond to
\(Q\)-graded commutative \(A\)-semialgebras.  The correspondence
commutes with arbitrary scalar extension.
\end{lemma}

\begin{proof}
The formal basis of \(A[Q]\) gives
\[
 B\ot_A A[Q]\cong\bigoplus_{q\in Q}B\,X^q.
\]
Write
\(\rho(b)=\sum_q b_q\ot X^q\).  The counit identity gives
\(b=\sum_qb_q\).  Coassociativity, compared coefficientwise in
\[
 B\ot_A A[Q]\ot_A A[Q]
 \cong\bigoplus_{p,q\in Q}B\,X^p\ot X^q,
\]
gives \(\rho(b_q)=b_q\ot X^q\).  Thus the \(b_q\) are homogeneous.
Their uniqueness follows from uniqueness of coefficients in the
formal monoid semiring.  Conversely, a grading defines
\(\rho(b_q)=b_q\ot X^q\), and the coaction identities are immediate.
The condition that \(\rho\) preserve multiplication is exactly
\(B_pB_q\subseteq B_{p+q}\).  Finally,
Lemma~\ref{lem:monoid-base-change} and distributivity of tensor
products over finite direct sums prove the base-change assertion.
\end{proof}

Let \(\operatorname{PicGr}_Q(A)\) be the groupoid of Picard-strong
\(Q\)-graded commutative \(A\)-semialgebras and graded
\(A\)-semialgebra isomorphisms.  Equivalently, an object consists of
invertible \(A\)-semimodules \(L_q\), an identification \(L_0\cong A\),
and coherent symmetric multiplication isomorphisms
\[
 L_p\ot_A L_q\longiso L_{p+q}.
\]
Thus it is a strong symmetric monoidal functor from the discrete
Picard group \(Q\) to the Picard groupoid of \(A\).

\subsection{Affine graded reconstruction}

\begin{theorem}[Graded reconstruction of diagonalizable torsors]
\label{thm:graded-torsor-equivalence}
Let \(Q\) be a finite abelian group of order \(d\), and let \(A\) be a
commutative semiring.  Taking the coordinate semialgebra and its
coaction induces an equivalence of groupoids
\begin{equation}
\label{eq:graded-torsor-equivalence}
 \Tors_{\fpqc}(\Spec A,D_A(Q))
 \longiso
 \operatorname{PicGr}_Q(A).
\end{equation}
More explicitly, every fpqc \(D_A(Q)\)-torsor is represented by an
affine morphism \(\Spec B\to\Spec A\), and
\[
 B=\bigoplus_{q\in Q}B_q
\]
is Picard-strong.  Consequently \(B\) is Zariski-locally free of rank
\(d\), and the torsor is finite locally free and faithfully flat.
Conversely, every Picard-strong grading satisfies the Hopf--Galois
identity
\begin{equation}
\label{eq:graded-hopf-galois}
 B\ot_A B
 \longiso
 B\ot_A A[Q].
\end{equation}
The equivalence and the decomposition commute with arbitrary base
change \(A\to A'\).
\end{theorem}

\begin{proof}
Let \(P\) be a \(D_A(Q)\)-torsor.  Choose an affine fpqc cover
\(\Spec A'\to\Spec A\) which trivializes \(P\).  Effective descent for
commutative algebras, Proposition~\ref{prop:descent-package}, produces
an \(A\)-semialgebra \(B\) representing \(P\).  The action gives an
\(A[Q]\)-coaction on \(B\), hence a grading
\[
 B=\bigoplus_{q\in Q}B_q
\]
by Lemma~\ref{lem:coaction-grading}.  An equivariant trivialization of
the torsor is a graded isomorphism
\[
 B\ot_A A'\longiso A'[Q].
\]
Taking homogeneous components gives
\[
 B_q\ot_A A'\longiso A'X^q
\]
for every \(q\).  Hence each \(B_q\) is fpqc-locally free of rank one.
The rank-one theorem of Borger and Jun
\cite[Theorem~1.1]{BorgerJun2025} makes \(B_q\) invertible and
Zariski-locally free of rank one.

The map \(A\to B_0\) and the multiplication maps
\[
 B_p\ot_A B_q\longrightarrow B_{p+q}
\]
become isomorphisms over \(A'\), because they do so in \(A'[Q]\).
They are therefore isomorphisms by the conservativity in
Proposition~\ref{prop:descent-package}.  Thus the grading is
Picard-strong.  Since \(Q\) is finite, a common Zariski refinement
trivializes all the \(B_q\); on that refinement \(B\) is free of rank
\(d\).  This proves the asserted finite local freeness without using
any higher-rank local-freeness descent statement.

Conversely, let \(B=\bigoplus_qB_q\) be Picard-strong.  It is finite
locally free of rank \(d\).  It is also faithfully flat: tensoring
with \(B\) is the finite direct sum of the equivalences
\(B_q\ot_A-\), so it preserves finite limits.  Moreover, for every
morphism \(f\) the map \(B\ot_A f\) is block diagonal for the
canonical \(Q\)-grading, and its degree-zero block is
\(B_0\ot_A f\cong f\).  If \(B\ot_A f\) is an isomorphism, uniqueness
of homogeneous decompositions makes its inverse block diagonal as
well; hence \(f\) is an isomorphism.  Thus \(B\) is flat and
conservative in the sense fixed in Section~\ref{sec:foundations}.
The coaction from Lemma~\ref{lem:coaction-grading} defines a
\(D_A(Q)\)-action.  The canonical map in
\eqref{eq:graded-hopf-galois} is
\begin{equation}
\label{eq:graded-canonical-map}
 \operatorname{can}(b\ot c_q)=bc_q\ot X^q.
\end{equation}
For fixed \(q\), its \(X^q\)-component is
\[
 B\ot_A B_q
 =
 \bigoplus_{p\in Q}B_p\ot_A B_q
 \longrightarrow
 \bigoplus_{p\in Q}B_{p+q}=B.
\]
This is an isomorphism because every map
\eqref{eq:strong-multiplication} is an isomorphism.  Summing over
\(q\) proves that \(\operatorname{can}\) is an isomorphism.  Together
with faithful flatness, this is the torsor identity.

The two constructions are inverse on objects and arrows because both
recover the same coaction.  Scalar extension preserves finite direct
sums, homogeneous components, invertibility, and the multiplication
maps, proving the last assertion.
\end{proof}

\begin{remark}
\label{rem:higher-rank-gap-resolved}
The proof does not assert that an arbitrary fpqc form of \(A^d\) is
Zariski-locally free.  It uses the \(D_A(Q)\)-coaction to split the
form canonically into \(d\) fpqc forms of \(A\), applies rank-one
descent to those pieces, and then takes their finite direct sum.  The
diagonalizable structure is therefore essential.
\end{remark}

\subsection{Global form and classifying geometry}

For a semiring scheme \(X\), let \(\operatorname{PicGr}_Q(X)\) denote
the groupoid of families of line bundles
\((\mathcal L_q)_{q\in Q}\), together with
\(\mathcal L_0\cong\OO_X\) and coherent symmetric isomorphisms
\[
 \mathcal L_p\ot\mathcal L_q\longiso\mathcal L_{p+q}.
\]
The following representability calculation is the semiring analogue
of the standard classifying-stack argument over rings
\cite{LaumonMoretBailly2000,Olsson2016}; its finite-local conclusion
comes from graded reconstruction rather than a general semiring
algebraicity theorem.

\begin{theorem}[Global graded descent]
\label{thm:global-graded-descent}
For every semiring scheme \(X\) and finite abelian group \(Q\), there
is a natural equivalence
\begin{equation}
\label{eq:global-graded-descent}
 \Tors_{\fpqc}(X,D_X(Q))
 \longiso
 \operatorname{PicGr}_Q(X).
\end{equation}
Every torsor in the left-hand groupoid is represented by a finite
locally free \(X\)-scheme of rank \(|Q|\).  If
\(\pi\colon P\to X\) is the corresponding torsor, then
\[
 \pi_*\OO_P=\bigoplus_{q\in Q}\mathcal L_q
\]
is its canonical Picard-strong grading.
\end{theorem}

\begin{proof}
On an affine open \(U=\Spec A\subseteq X\), affine fpqc descent
represents the restriction of a torsor, and
Theorem~\ref{thm:graded-torsor-equivalence} gives its canonical
grading.  The construction is functorial under localization, so the
affine coordinate algebras and their homogeneous pieces glue.
Conversely,
\[
 \mathcal B=\bigoplus_{q\in Q}\mathcal L_q
\]
is a quasi-coherent Picard-strong graded semialgebra.  Its relative
spectrum is a torsor by the affine theorem.  These constructions are
inverse Zariski-locally and hence globally.  Finite local freeness of
rank \(|Q|\) follows because the finitely many line bundles
\(\mathcal L_q\) are simultaneously trivial on a common Zariski
refinement.
\end{proof}

\begin{corollary}[Finite diagonalizable classifying stack]
\label{cor:general-classifying-stack}
Let \(Q\) be finite of order \(d\).  The trivial torsor defines a
representable finite locally free fpqc-surjective morphism
\[
 X\longrightarrow\mathscr B_XD_X(Q)
\]
of rank \(d\).  The diagonal of \(\mathscr B_XD_X(Q)\) is
representable and finite locally free of rank \(d\).  Both assertions
commute with arbitrary base change.
\end{corollary}

\begin{proof}
After a map \(T\to\mathscr B_XD_X(Q)\) corresponding to a torsor
\(P\to T\), the first morphism becomes the sheaf of trivializations
of \(P\), which is \(P\) itself.  It is finite locally free of rank
\(d\) by Theorem~\ref{thm:global-graded-descent}.  The fiber of the
diagonal over two torsors \(P_1,P_2\) is
\(\Isom_{D(Q)}(P_1,P_2)\), again a \(D_T(Q)\)-torsor.  Apply the same
theorem.  Base-change compatibility follows from
Theorem~\ref{thm:graded-torsor-equivalence}.
\end{proof}

\begin{corollary}[Homogeneous Picard invariant]
\label{cor:homogeneous-picard-invariant}
A \(D_X(Q)\)-torsor \(P\) determines a homomorphism
\[
 \omega_P\colon Q\longrightarrow\Pic(X),
 \qquad
 q\longmapsto[\mathcal L_q].
\]
Its image is killed by \(|Q|\).  If \(\Pic(X)\) has no
\(|Q|\)-torsion, then all homogeneous line bundles are trivial,
although the torsor can remain nontrivial through its multiplication
constants.  The affine Boolean torsors of
Theorem~\ref{thm:affine-boolean-kummer-torsor} exhibit this latter
possibility explicitly.
\end{corollary}

\begin{proof}
The strong multiplication isomorphisms give
\([\mathcal L_{p+q}]=[\mathcal L_p]+[\mathcal L_q]\).
Since \(|Q|q=0\), the image of every \(q\) is killed by \(|Q|\).
\end{proof}

\begin{corollary}[Cyclic reconstruction]
\label{cor:cyclic-graded-reconstruction}
For \(Q=C_m\), Picard-strong \(C_m\)-gradings are naturally
equivalent to pairs
\[
 (L,\vartheta),\qquad
 L\ \text{a line bundle},\qquad
 \vartheta\colon L^{\ot m}\longiso\OO_X.
\]
The corresponding graded algebra is
\[
 \OO_X\oplus L\oplus\cdots\oplus L^{\ot(m-1)},
\]
with multiplication wrapped by \(\vartheta\).
\end{corollary}

\begin{proof}
From a strong grading take \(L=\mathcal L_{\bar1}\).  Iterated
multiplication identifies \(\mathcal L_{\bar r}\) with \(L^{\ot r}\)
for \(0\leq r<m\), and the final multiplication gives
\(L^{\ot m}\cong\mathcal L_{\bar0}\cong\OO_X\).  Associativity makes
these identifications independent of parenthesization.  Conversely,
the displayed direct sum becomes strongly graded by using ordinary
tensor multiplication when \(i+j<m\) and applying \(\vartheta\) once
when \(i+j\geq m\).  The two constructions are inverse.
\end{proof}

\section{Matrix Kummer stacks and exact sequences}
\label{sec:matrix-kummer}

The power morphism is the one-dimensional instance of a lattice
isogeny.  We now retain the entire lattice presentation and obtain a
multi-parameter Kummer theory.  This adds a second description to the
graded model of Section~\ref{sec:graded-torsors}: a finite
diagonalizable torsor may be encoded either by all of its homogeneous
line bundles or by a finite presentation of its character group.

\subsection{Matrix notation and torus kernels}

Fix \(r\geq1\) and an integral matrix
\[
 M=(m_{ij})\in\operatorname{Mat}_{r\times r}(\mathbb Z),
 \qquad \det M\neq0.
\]
For an abelian group \(G\), written multiplicatively, set
\begin{equation}
\label{eq:matrix-power-group}
 \Phi_M(g_1,\ldots,g_r)_j
 =
 \prod_{i=1}^r g_i^{m_{ij}}
 \qquad(1\leq j\leq r).
\end{equation}
For a tuple of line bundles
\(\mathbf L=(L_1,\ldots,L_r)\), define
\begin{equation}
\label{eq:matrix-power-lines}
 \Phi_M(\mathbf L)_j
 =
 \bigotimes_{i=1}^r L_i^{\ot m_{ij}},
\end{equation}
where negative tensor powers mean powers of the dual.  The same
formula on points defines a homomorphism of split tori
\[
 \varphi_M\colon\Gm^r\longrightarrow\Gm^r.
\]
Its pullback on character lattices is \(M\).  Put
\begin{equation}
\label{eq:matrix-kernel}
 Q_M=\coker(M\colon\mathbb Z^r\to\mathbb Z^r),
 \qquad
 K_M=D(Q_M).
\end{equation}
Thus \(|Q_M|=|\det M|\).

\begin{proposition}[Finite torus isogeny]
\label{prop:matrix-finite-torsor}
The sequence
\begin{equation}
\label{eq:matrix-sheaf-sequence}
 1\longrightarrow K_M\longrightarrow\Gm^r
 \xrightarrow{\ \varphi_M\ }\Gm^r\longrightarrow1
\end{equation}
is exact on the fpqc site of semiring schemes.  The middle morphism is
finite locally free and faithfully flat of rank \(|\det M|\), and is
a \(K_M\)-torsor.  These assertions commute with arbitrary change of
the base semiring.
\end{proposition}

\begin{proof}
Apply Theorem~\ref{thm:finite-index-torsor} to the injective character
map \(M\colon\mathbb Z^r\to\mathbb Z^r\).  The kernel is
\(D(Q_M)\), the torsor identity gives exactness at the first two
terms, and faithful flatness makes \(\varphi_M\) an epimorphism of
fpqc sheaves.
\end{proof}

We use the term \emph{finite torus isogeny} only for a homomorphism
having the properties in Proposition~\ref{prop:matrix-finite-torsor};
no independent notion of dimension for semiring group schemes is
being invoked.

\subsection{The homotopy fiber of the matrix power functor}

Define \(\Kum_M(X)\) to be the following symmetric monoidal groupoid.
An object is a tuple
\begin{equation}
\label{eq:matrix-kummer-object}
 (\mathbf L,\boldsymbol\vartheta)
 =
 \left(
  (L_1,\ldots,L_r),
  (\vartheta_1,\ldots,\vartheta_r)
 \right),
\end{equation}
where the \(L_i\) are line bundles and
\begin{equation}
\label{eq:matrix-trivializations}
 \vartheta_j\colon
 \Phi_M(\mathbf L)_j\longiso\OO_X
 \qquad(1\leq j\leq r).
\end{equation}
A morphism
\((\mathbf L,\boldsymbol\vartheta)\to
(\mathbf L',\boldsymbol\vartheta')\)
is a tuple of line-bundle isomorphisms
\(\mathbf f=(f_i)\) satisfying
\[
 \vartheta_j
 =
 \vartheta'_j\circ\Phi_M(\mathbf f)_j
 \qquad(1\leq j\leq r).
\]
Tensor product is componentwise, with the symmetry constraints used
to collect equal line-bundle factors.  In other words,
\(\Kum_M(X)\) is the homotopy fiber over the trivial tuple of
\[
 \Phi_M\colon\operatorname{Picard}(X)^r
 \longrightarrow\operatorname{Picard}(X)^r.
\]

\begin{theorem}[Matrix Kummer equivalence]
\label{thm:matrix-kummer-equivalence}
For every semiring scheme \(X\), there are natural symmetric
monoidal equivalences
\begin{equation}
\label{eq:three-kummer-models}
 \Kum_M(X)
 \longiso
 \Tors_{\fpqc}(X,K_M)
 \longiso
 \operatorname{PicGr}_{Q_M}(X).
\end{equation}
Under the second equivalence, every matrix Kummer object determines a
finite locally free \(K_M\)-torsor of rank \(|\det M|\), whose
coordinate algebra is the Picard-strong \(Q_M\)-graded algebra of
Theorem~\ref{thm:global-graded-descent}.
\end{theorem}

\begin{proof}
We first prove the left-hand equivalence.  More generally, consider
an exact sequence of abelian fpqc sheaves
\[
 1\longrightarrow K\longrightarrow T
 \xrightarrow{\varphi}T'\longrightarrow1.
\]
A \(K\)-torsor \(P\) extends to a \(T\)-torsor
\[
 E=P\times^K T.
\]
Its induced \(T'\)-torsor is canonically trivial, because
\[
 E\times^T T'
 \cong
 P\times^K T'
 \cong
 (P/K)\times T'
 \cong X\times T';
\]
the action of \(K\) on \(T'\) is trivial.  Conversely, if a
\(T\)-torsor \(E\) is equipped with a trivialization of its induced
\(T'\)-torsor, the inverse image in \(E\) of the identity section of
\(T'\) is a \(K\)-torsor.  On a cover trivializing \(E\), these are
the evident inverse constructions; descent makes them inverse
globally.

Apply this observation to
\eqref{eq:matrix-sheaf-sequence}.  A \(\Gm^r\)-torsor is the product
of the coframe torsors of a tuple
\(\mathbf L=(L_1,\ldots,L_r)\).  Under extension by
\(\varphi_M\), its \(j\)th coframe torsor is the coframe torsor of
\(\Phi_M(\mathbf L)_j\).  A trivialization of the induced
\(\Gm^r\)-torsor is therefore exactly the tuple
\(\boldsymbol\vartheta\) in
\eqref{eq:matrix-trivializations}.  This proves
\(\Kum_M(X)\simeq\Tors_{\fpqc}(X,K_M)\), including the stated
morphisms and tensor products.

The second equivalence and finite local freeness follow from
Theorem~\ref{thm:global-graded-descent}, because
\(K_M=D(Q_M)\).
\end{proof}

\begin{remark}
\label{rem:matrix-conventions}
The use of coframes \(L_i\to\OO_X\) fixes all signs in
\eqref{eq:matrix-power-lines}.  With the alternative convention of
frames \(\OO_X\to L_i\), every character is inverted.  The resulting
groupoids are canonically equivalent by dualizing all line bundles,
but mixing the two conventions would incorrectly replace \(M\) by
\(-M\).
\end{remark}

\subsection{The matrix Kummer exact sequence}

For an abelian group \(G\), write
\[
 \coker\Phi_M(G)
 =
 G^r/\Phi_M(G^r),
 \qquad
 \ker\Phi_M(G)
 =
 \ker\bigl(\Phi_M\colon G^r\to G^r\bigr).
\]
On \(\Pic(X)\) we use additive notation, so
\[
 \Phi_M(\ell_1,\ldots,\ell_r)_j
 =
 \sum_{i=1}^r m_{ij}\ell_i.
\]

\begin{theorem}[Matrix Kummer exact sequence]
\label{thm:matrix-kummer-exact}
Let \(U_X=\Gamma(X,\OO_X^\times)\).  There is a natural short exact
sequence of abelian groups
\begin{equation}
\label{eq:matrix-kummer-exact}
 0\longrightarrow
 \coker\Phi_M(U_X)
 \xrightarrow{\ \delta_{M,X}\ }
 H^1_{\fpqc}(X,K_M)
 \xrightarrow{\ \pi_{M,X}\ }
 \ker\Phi_M(\Pic(X))
 \longrightarrow0.
\end{equation}
The map \(\pi_{M,X}\) forgets the trivializations in a matrix Kummer
object.  For \(\mathbf u\in U_X^r\), the boundary
\(\delta_{M,X}(\mathbf u)\) is the pullback of the universal
\(K_M\)-torsor
\(\varphi_M\colon\Gm^r\to\Gm^r\)
along the point \(\mathbf u\colon X\to\Gm^r\).
\end{theorem}

\begin{proof}
Use the first equivalence in
Theorem~\ref{thm:matrix-kummer-equivalence}.  The underlying tuple
\(\mathbf L\) of an object
\((\mathbf L,\boldsymbol\vartheta)\) satisfies
\(\Phi_M([\mathbf L])=0\), which defines \(\pi_{M,X}\).  Conversely,
if \(\Phi_M([\mathbf L])=0\), each line bundle
\(\Phi_M(\mathbf L)_j\) is trivial; choosing trivializations produces
an object of \(\Kum_M(X)\).  Hence \(\pi_{M,X}\) is surjective.

An object lies in its kernel exactly when all the \(L_i\) are
trivial.  After choosing trivializations of them, the tuple
\(\boldsymbol\vartheta\) is multiplication by a vector
\(\mathbf u\in U_X^r\).  Changing the chosen trivializations by
\(\mathbf a\in U_X^r\) multiplies \(\mathbf u\) by an element of
\(\Phi_M(U_X^r)\).  Thus the kernel is
\(\coker\Phi_M(U_X)\).  For the trivial tuple of line
bundles, the inverse-image construction in the proof of
Theorem~\ref{thm:matrix-kummer-equivalence} is precisely the pullback
of \(\varphi_M\) along \(\mathbf u\).  This identifies the first map
with the stated boundary.
\end{proof}

\begin{corollary}[Matrix classifying stack]
\label{cor:matrix-classifying-stack}
The stack \(T\mapsto\Kum_M(T)\) is naturally
\(\mathscr B_XK_M\).  It admits a representable finite locally free
fpqc-surjective atlas
\[
 X\longrightarrow\mathscr B_XK_M
\]
of rank \(|\det M|\), and its diagonal is representable and finite
locally free of the same rank.  The equivalence, atlas, diagonal, and
exact sequence commute with arbitrary base change \(X'\to X\).
\end{corollary}

\begin{proof}
The stack identification is
Theorem~\ref{thm:matrix-kummer-equivalence}; the geometric assertions
are Corollary~\ref{cor:general-classifying-stack} for \(Q=Q_M\).
Pullback commutes with tensor products, duals, trivializations,
gradings, and relative spectra, proving base-change compatibility.
\end{proof}

\subsection{Smith normal form}

\begin{theorem}[Invariant-factor decomposition]
\label{thm:smith-kummer}
Let
\[
 UMV=\operatorname{diag}(d_1,\ldots,d_r)
\]
be a Smith normal form, where \(U,V\in\operatorname{GL}_r(\mathbb Z)\),
\(d_i\geq1\), and \(d_i\mid d_{i+1}\).  The induced changes of torus
coordinates give noncanonical isomorphisms
\begin{equation}
\label{eq:smith-kernel}
 K_M\cong\prod_{i=1}^r\mu_{d_i},
\qquad
 \mathscr B_XK_M\cong
 \prod_{i=1}^r\mathscr B_X\mu_{d_i}.
\end{equation}
Consequently,
\begin{equation}
\label{eq:smith-H1}
 H^1_{\fpqc}(X,K_M)
 \cong
 \prod_{i=1}^rH^1_{\fpqc}(X,\mu_{d_i}).
\end{equation}
The finite locally free rank of the universal matrix Kummer cover is
\[
 \prod_{i=1}^r d_i=|\det M|.
\]
\end{theorem}

\begin{proof}
Unimodular matrices induce automorphisms of \(\Gm^r\).  Changing the
source and target character bases by \(U\) and \(V\) identifies the
kernel of \(\varphi_M\) with the kernel of the diagonal morphism
\[
 (x_1,\ldots,x_r)\longmapsto
 (x_1^{d_1},\ldots,x_r^{d_r}).
\]
That kernel is \(\prod_i\mu_{d_i}\).  Products of commutative torsors
give \eqref{eq:smith-H1}, and the rank formula follows either from
Proposition~\ref{prop:matrix-finite-torsor} or from the product of the
cyclic ranks.
\end{proof}

\begin{remark}
The decomposition depends on Smith bases.  In contrast,
\(\operatorname{PicGr}_{Q_M}(X)\), the stack \(\mathscr B_XK_M\),
and \eqref{eq:matrix-kummer-exact} are intrinsic to the lattice
presentation.  The matrix theorem therefore retains the functorial
relation between the original torus coordinates, global units, and
line bundles.
\end{remark}

\section{An explicit affine torsor in characteristic one}
\label{sec:affine-boolean-torsor}

The preceding constructions are not confined to projective examples or
to torsors arising from nontrivial line bundles.  We now give an affine
example in which the unit term of the Kummer sequence is nonzero and the
coordinate semialgebra is completely explicit.  It also shows that a
group scheme can have only the identity as a base-valued point while
possessing a nontrivial torsor.

Fix \(m\geq2\), and put
\[
 A=\mathbb B[u,u^{-1}]=\mathbb B[\mathbb Z],
 \qquad
 B=\mathbb B[v,v^{-1}]=\mathbb B[\mathbb Z].
\]
Regard \(B\) as an \(A\)-semialgebra through
\begin{equation}
\label{eq:boolean-affine-power-map}
 \iota_m\colon A\longrightarrow B,
 \qquad u\longmapsto v^m.
\end{equation}

\begin{lemma}[Units of the Boolean Laurent semiring]
\label{lem:boolean-laurent-units}
Every unit of \(A\) is a Laurent monomial.  More precisely,
\[
 A^\times=\{u^n:n\in\mathbb Z\}\cong\mathbb Z.
\]
Consequently,
\[
 A^\times/(A^\times)^m\cong\mathbb Z/m\mathbb Z,
\]
with the class of \(u\) as a generator.
\end{lemma}

\begin{proof}
For a nonzero Boolean Laurent polynomial \(f\), let
\(\operatorname{supp}(f)\subset\mathbb Z\) be its finite nonempty
support.  Boolean addition identifies \(f\) with its support, and
multiplication gives
\[
 \operatorname{supp}(fg)
 =\operatorname{supp}(f)+\operatorname{supp}(g).
\]
If \(fg=1\), the Minkowski sum on the right is the singleton
\(\{0\}\).  Comparing the least and greatest exponents shows that both
supports have width zero; hence each is a singleton.  Thus
\(f=u^n\) and \(g=u^{-n}\) for a unique \(n\).  The quotient statement
is now the cokernel of multiplication by \(m\) on \(\mathbb Z\).
\end{proof}

\begin{theorem}[Affine Boolean Kummer torsor]
\label{thm:affine-boolean-kummer-torsor}
Let
\[
 X=\operatorname{Spec}A,
 \qquad P_m=\operatorname{Spec}B.
\]
Then \eqref{eq:boolean-affine-power-map} defines a finite free
\(\mu_m\)-torsor
\[
 \pi_m\colon P_m\longrightarrow X
\]
of rank \(m\).  As an \(A\)-semimodule,
\[
 B=A\oplus Av\oplus\cdots\oplus Av^{m-1},
 \qquad v^m=u,
\]
and this is its Picard-strong \(C_m\)-grading.  The torsor has no
section over \(X\), and hence is nontrivial.  Its cohomology class is
the Kummer boundary \(\delta_m(u)\).
\end{theorem}

\begin{proof}
The map \(\iota_m\) is the coordinate map of the power morphism
\[
 [m]\colon\mathbf G_{\mathrm m,\mathbb B}\longrightarrow
 \mathbf G_{\mathrm m,\mathbb B}.
\]
Theorem~\ref{thm:power-torsor}
and Proposition~\ref{prop:coset-decomposition} therefore give the
torsor assertion and the displayed free decomposition.  Equivalently,
it is the root semialgebra
\[
 A\langle u^{1/m}\rangle
 =A[z,z^{-1}]/(z^m\sim u),
\]
with \(z\) identified with \(v\).

A section of \(\pi_m\) would be an \(A\)-algebra map \(B\to A\), and
the image \(w\) of \(v\) would satisfy \(w^m=u\).  By
Lemma~\ref{lem:boolean-laurent-units}, \(w=u^n\) for some
\(n\in\mathbb Z\), which would force \(mn=1\), impossible for
\(m\geq2\).  A torsor is trivial if and only if it has a global
section, so \(\pi_m\) is nontrivial.  Its identification with
\(\delta_m(u)\) follows from the fiber description in
Theorem~\ref{thm:matrix-kummer-exact}, specialized to the matrix
\((m)\).
\end{proof}

\begin{corollary}[A visible nonzero subgroup of \(H^1\)]
\label{cor:boolean-affine-H1}
The cyclic Kummer sequence induces a canonical injection
\[
 \mathbb Z/m\mathbb Z
 \cong\Ext^1_{\mathbb Z}(C_m,A^\times)
 \cong A^\times/(A^\times)^m
 \lhook\joinrel\longrightarrow
 H^1_{\fpqc}(X,\mu_m).
\]
For \(a\in\mathbb Z\), the image of \(\bar a\) is represented by the
finite free root semialgebra
\[
 A\langle u^{a/m}\rangle
 =\bigoplus_{i=0}^{m-1}Az^i,
 \qquad z^m=u^a,
\]
with the wrapped multiplication of
\eqref{eq:wrapped-multiplication}.  These torsors are isomorphic
exactly when their exponents are congruent modulo \(m\).
\end{corollary}

\begin{proof}
The beginning of the long exact sequence associated with
\[
 1\longrightarrow\mu_m\longrightarrow\Gm
 \xrightarrow{[m]}\Gm\longrightarrow1
\]
injects \(A^\times/(A^\times)^m\) into
\(H^1_{\fpqc}(X,\mu_m)\).  Lemma~\ref{lem:boolean-laurent-units}
identifies this quotient with \(\mathbb Z/m\mathbb Z\), and
Proposition~\ref{prop:root-algebra} gives the displayed representatives.
Injectivity gives the final isomorphism criterion.
\end{proof}

\begin{remark}[Characteristic-one content]
\label{rem:boolean-affine-shadow}
The example is affine, non-ring-theoretic, and independent of any
Picard-group computation.  Since \(A\) is a nonzero idempotent
semiring, its additive group completion is the zero ring.  Thus the
nonzero class in Corollary~\ref{cor:boolean-affine-H1} has only a
degenerate classical shadow, just as the Boolean root gerbe in
Corollary~\ref{cor:boolean-cohomological-jump}.
\end{remark}

\section{Derived character cohomology}
\label{sec:derived}

The graded theorem identifies the geometry of degree-one
diagonalizable classes.  We now record the cohomological consequences
of the finite torus resolution.  The coefficient objects are ordinary
sheaves of abelian groups: no derived category of semimodules is used,
and the results below reduce diagonalizable cohomology to
\(\Gm\)-cohomology rather than evaluating the latter.

\subsection{A derived character formula}

For an abelian fpqc sheaf \(\mathcal F\) on a semiring scheme \(X\),
write
\[
 H^i_{\fpqc}(X,\mathcal F)
 =
 R^i\Gamma(X_{\fpqc},\mathcal F).
\]
The category of abelian sheaves on a site has enough injectives, so
these groups and their long exact sequences have their usual
meaning~\cite{StacksCohomologySites}.  In degree one this agrees with
the torsor group used above; in degree two it classifies gerbes
banded by \(\mathcal F\) when \(\mathcal F\) is abelian
\cite[Chapter~IV]{Giraud1971}.

\begin{lemma}[Character resolution]
\label{lem:character-resolution}
Let \(Q\) be a finite abelian group.  There is a free resolution
\begin{equation}
\label{eq:character-resolution}
 0\longrightarrow\mathbb Z^r
 \xrightarrow{\ M\ }\mathbb Z^r
 \longrightarrow Q\longrightarrow0
\end{equation}
with \(\det M\neq0\).  For every abelian group \(C\), precomposition
with \(M\) gives an exact sequence
\begin{equation}
\label{eq:hom-ext-resolution}
 0\longrightarrow\Hom_{\mathbb Z}(Q,C)
 \longrightarrow C^r
 \xrightarrow{\ \Phi_M\ } C^r
 \longrightarrow\Ext^1_{\mathbb Z}(Q,C)
 \longrightarrow0.
\end{equation}
In particular,
\[
 \ker\Phi_M(C)=\Hom_{\mathbb Z}(Q,C),
 \qquad
 \coker\Phi_M(C)=\Ext^1_{\mathbb Z}(Q,C).
\]
\end{lemma}

\begin{proof}
Choose invariant factors
\[
 Q\cong\bigoplus_{j=1}^r\mathbb Z/d_j\mathbb Z
\]
and take \(M=\operatorname{diag}(d_1,\ldots,d_r)\).  More generally,
any full-rank square presentation can be used.  Apply
\(\Hom_{\mathbb Z}(-,C)\) to
\eqref{eq:character-resolution}.  The two free groups have vanishing
\(\Ext^1\), which gives \eqref{eq:hom-ext-resolution}.  In chosen
bases the middle map is
\[
 (c_i)_i\longmapsto
 \left(\sum_i m_{ij}c_i\right)_j;
\]
this is the additive form of \(\Phi_M\).
\end{proof}

\begin{theorem}[Relative derived character formula]
\label{thm:relative-derived-character-formula}
Let \(f\colon X\to Y\) be a morphism of semiring schemes and let
\(Q\) be a finite abelian group.  If \(Q_Y\) denotes the associated
constant abelian sheaf on \(Y_{\fpqc}\), then there is a natural
isomorphism in the derived category of abelian fpqc sheaves on \(Y\):
\begin{equation}
\label{eq:relative-derived-character-formula}
 Rf_*D_X(Q)
 \longiso
 R\mathcal{H}\!om_{\mathbb Z,Y}
 \bigl(Q_Y,Rf_*\Gm{}_{X}\bigr).
\end{equation}
It is contravariantly functorial in \(Q\), pseudofunctorial in \(f\)
under the canonical composition isomorphisms for derived direct image,
and independent of a presentation of \(Q\).
\end{theorem}

\begin{proof}
Let
\(f_{\fpqc}\colon X_{\fpqc}\to Y_{\fpqc}\) be the geometric
morphism of fpqc topoi induced by \(f\); all direct images in this
proof refer to this morphism.  Choose
\eqref{eq:character-resolution}.  Applying \(D_X(-)\) gives
the sequence of abelian fpqc sheaves
\begin{equation}
\label{eq:character-sheaf-resolution}
 1\longrightarrow D_X(Q)\longrightarrow(\Gm{}_{X})^r
 \xrightarrow{\ \varphi_M\ }(\Gm{}_{X})^r\longrightarrow1.
\end{equation}
It is exact by Proposition~\ref{prop:matrix-finite-torsor}; in
particular, the last arrow is an epimorphism of fpqc sheaves.  Thus,
with the source torus in degree zero and the target torus in degree
one, the associated two-term complex is quasi-isomorphic to its
kernel:
\[
 D_X(Q)\longiso
 [\,(\Gm{}_{X})^r\xrightarrow{\varphi_M}(\Gm{}_{X})^r\,]
\]
in the derived category of abelian fpqc sheaves on \(X\).

Apply \(Rf_*\) to the resulting distinguished triangle.  Since
direct image commutes with finite products, this identifies
\(Rf_*D_X(Q)\) with the homotopy fiber of
\[
 (Rf_*\Gm{}_X)^r\xrightarrow{\varphi_M}(Rf_*\Gm{}_X)^r.
\]

For comparison on the right-hand side, let
\[
 P^\bullet=
 [\,\mathbb Z_Y^r\xrightarrow{M}\mathbb Z_Y^r\,]
\]
with terms in degrees \(-1\) and \(0\).  The constant-sheaf functor is
exact on abelian groups
\cite[Lemma~18.42.1]{StacksConstantSheaves}; hence sheafification of
\eqref{eq:character-resolution} is exact and
\(P^\bullet\to Q_Y\) is a finite free resolution.  Moreover,
for every abelian sheaf \(\mathcal F\),
\[
 \mathcal{H}\!om_{\mathbb Z,Y}(\mathbb Z_Y^r,\mathcal F)
 \cong\mathcal F^r,
\]
and this functor is exact.  We do not need to assert that
\(\mathbb Z_Y\) is projective in the abelian category of sheaves: the
displayed exactness is precisely what permits the bounded complex
\(P^\bullet\) to compute derived internal Hom.  Hence
\[
 R\mathcal{H}\!om_{\mathbb Z,Y}
 \bigl(Q_Y,Rf_*\Gm{}_{X}\bigr)
 \longiso
 \mathcal{H}\!om_{\mathbb Z,Y}
 \bigl(P^\bullet,Rf_*\Gm{}_{X}\bigr).
\]
The last complex is the same homotopy fiber: precomposition with
\(M\) is the monomial map \(\varphi_M\) after identifying internal
Hom from \(\mathbb Z_Y^r\) with an \(r\)-fold product.  This proves
\eqref{eq:relative-derived-character-formula}, including its degree
conventions.

A homomorphism of finite abelian groups lifts, before passage to
constant sheaves, to a chain map between free abelian resolutions, and
two lifts are chain homotopic.  Sheafification preserves the resulting
maps and homotopies.  The comparison theorem for projective
resolutions therefore proves functoriality in \(Q\) and independence
of all choices; see, for example, \cite[Chapters~2 and~5]{Weibel1994}.
The construction uses only the
canonical geometric morphism of fpqc topoi and the natural unit maps
for direct image.  It is consequently compatible with commutative
triangles of semiring schemes and, up to the canonical derived
composition isomorphisms, with composition of direct-image functors.
This is the asserted pseudofunctoriality in \(f\).
\end{proof}

\begin{corollary}[Derived character formula]
\label{thm:derived-character-formula}
Let \(X\) be a semiring scheme and \(Q\) a finite abelian group.
There is a natural isomorphism in the derived category of abelian
groups
\begin{equation}
\label{eq:derived-character-formula}
 R\Gamma_{\fpqc}\!\left(X,D_X(Q)\right)
 \longiso
 R\Hom_{\mathbb Z}\!\left(
 Q,R\Gamma_{\fpqc}(X,\Gm)
 \right).
\end{equation}
It is contravariantly functorial in \(Q\), functorial under pullback
in \(X\), and independent of a presentation of \(Q\).
\end{corollary}

\begin{proof}
Apply Theorem~\ref{thm:relative-derived-character-formula} to the
direct-image functor from \(X_{\fpqc}\) to the punctual topos.
\end{proof}

\subsection{The universal-coefficient filtration}

\begin{theorem}[Relative universal-coefficient filtration]
\label{thm:relative-universal-coefficient}
Under the hypotheses of
Theorem~\ref{thm:relative-derived-character-formula}, for every
\(i\geq1\) there is a natural short exact sequence of abelian sheaves
on \(Y_{\fpqc}\):
\begin{equation}
\label{eq:relative-universal-coefficient}
 0\longrightarrow
 \mathcal{E}\!xt^1_{\mathbb Z,Y}
 \bigl(Q_Y,R^{i-1}f_*\Gm{}_{X}\bigr)
 \longrightarrow
 R^if_*D_X(Q)
\end{equation}
\[
 \longrightarrow
 \mathcal{H}\!om_{\mathbb Z,Y}
 \bigl(Q_Y,R^if_*\Gm{}_{X}\bigr)
 \longrightarrow0.
\]
This is the two-step filtration induced by
\eqref{eq:relative-derived-character-formula}; it need not split
naturally.
\end{theorem}

\begin{proof}
Take cohomology sheaves in
\eqref{eq:relative-derived-character-formula}.  Since \(Q_Y\) has the
length-one free resolution
\eqref{eq:character-resolution}, the hyper-Ext spectral sequence has
only the columns \(0\) and \(1\).  Its edge filtration is exactly
\eqref{eq:relative-universal-coefficient}.  Equivalently, apply
\(Rf_*\) to \eqref{eq:character-sheaf-resolution}, take cohomology
sheaves, and identify the resulting kernels and cokernels by the
sheafified form of Lemma~\ref{lem:character-resolution}.
\end{proof}

\begin{theorem}[Universal coefficient theorem]
\label{thm:universal-coefficient}
For every \(i\geq1\) there is a natural short exact sequence
\begin{equation}
\label{eq:universal-coefficient}
\begin{aligned}
 0&\longrightarrow
 \Ext^1_{\mathbb Z}\!\left(Q,H^{i-1}_{\fpqc}(X,\Gm)\right)
 \xrightarrow{\ \alpha_i\ }
 H^i_{\fpqc}\!\left(X,D_X(Q)\right)\\
 &\xrightarrow{\ \beta_i\ }
 \Hom_{\mathbb Z}\!\left(Q,H^i_{\fpqc}(X,\Gm)\right)
 \longrightarrow0.
\end{aligned}
\end{equation}
The sequence is natural in \(X\) and \(Q\), but it need not split
naturally.  The right-hand map sends a class \(\xi\) to the
homomorphism
\[
 q\longmapsto q_*(\xi),
\]
where \(q\colon D_X(Q)\to\Gm\) is the character indexed by
\(q\in Q\).
\end{theorem}

\begin{proof}
Specialize
Theorem~\ref{thm:relative-universal-coefficient} to the direct image
to the punctual topos.  The description of \(\beta_i\) follows by
pushing a class forward along all characters.
\end{proof}

\begin{corollary}[Invariant degree-one sequence]
\label{cor:invariant-degree-one}
Let \(U_X=\Gamma(X,\OO_X^\times)\).  There is a natural exact
sequence
\begin{equation}
\label{eq:invariant-degree-one}
 0\longrightarrow
 \Ext^1_{\mathbb Z}(Q,U_X)
 \longrightarrow
 H^1_{\fpqc}(X,D_X(Q))
 \longrightarrow
 \Hom_{\mathbb Z}(Q,\Pic(X))
 \longrightarrow0.
\end{equation}
For a presentation \(Q=\coker M\), it is exactly the
matrix Kummer sequence
\eqref{eq:matrix-kummer-exact}.
\end{corollary}

\begin{proof}
Use \(H^0_{\fpqc}(X,\Gm)=U_X\),
\(H^1_{\fpqc}(X,\Gm)=\Pic(X)\), and
Lemma~\ref{lem:character-resolution}.
\end{proof}

\begin{theorem}[Realization and ambiguity of homogeneous Picard data]
\label{thm:realization-ambiguity}
Let \(Q\) be finite.  Every homomorphism
\[
 \omega\colon Q\longrightarrow\Pic(X)
\]
is realized by a Picard-strong \(Q\)-graded semialgebra and hence by a
finite locally free \(D_X(Q)\)-torsor.  For a fixed \(\omega\), the set
of isomorphism classes of realizations is a principal homogeneous
space under
\[
 \Ext^1_{\mathbb Z}(Q,U_X),
 \qquad U_X=\Gamma(X,\OO_X^\times).
\]
The action changes the coherent homogeneous multiplication data while
leaving the classes of the homogeneous line bundles fixed.
The automorphism group of every realization is canonically
\[
 D_X(Q)(X)=\Hom_{\mathbb Z}(Q,U_X).
\]
In the graded model, a homomorphism \(\chi\colon Q\to U_X\) acts on
the summand of degree \(q\) by multiplication with \(\chi(q)\).

Consequently:
\begin{enumerate}
\item if there are no nonzero homomorphisms from \(Q\) to
\(\Pic(X)\), then all torsor classes come from unit-extension data and
\[
 H^1_{\fpqc}(X,D_X(Q))
 \cong\Ext^1_{\mathbb Z}(Q,U_X);
\]
\item if \(\Ext^1_{\mathbb Z}(Q,U_X)=0\), then the homogeneous Picard
invariant classifies torsors uniquely up to isomorphism.
\end{enumerate}
\end{theorem}

\begin{proof}
Surjectivity and the description of every fiber follow from the short
exact sequence \eqref{eq:invariant-degree-one}.  To make realization
explicit, choose a length-one free presentation with equal-rank terms,
\[
 0\longrightarrow\mathbb Z^r\xrightarrow{M}\mathbb Z^r
 \longrightarrow Q\longrightarrow0.
\]
Compose the basis elements of the right-hand free group with \(\omega\)
and choose line bundles \(L_1,\ldots,L_r\) representing the resulting
classes.  The relations encoded by \(M\) say precisely that every
\(\Phi_M(\mathbf L)_j\) has trivial Picard class.  Choosing
trivializations produces a matrix Kummer object
\((\mathbf L,\boldsymbol\vartheta)\).  Theorem
\ref{thm:matrix-kummer-equivalence} turns it into a finite locally free
torsor and its Picard-strong grading.  Two choices over the same
\(\omega\) differ by the kernel of the forgetful map, which
Theorem~\ref{thm:matrix-kummer-exact} and
Lemma~\ref{lem:character-resolution} identify with
\(\Ext^1_{\mathbb Z}(Q,U_X)\).  The final assertions are the two
extreme cases of the same exact sequence.  Finally, automorphisms of a
torsor under the commutative group \(D_X(Q)\) are translations by
global sections of that group.  The displayed description follows
from \(D_X(Q)(X)=\Hom_{\mathbb Z}(Q,U_X)\), and the formula on
homogeneous summands is the corresponding graded action.
\end{proof}

\begin{corollary}[Exponent bound]
\label{cor:derived-exponent}
If \(e\) is the exponent of \(Q\), then
\[
 e\,H^i_{\fpqc}(X,D_X(Q))=0
\qquad(i\geq0).
\]
\end{corollary}

\begin{proof}
Multiplication by \(e\) on the sheaf \(D_X(Q)\) is zero, because it
is dual to multiplication by \(e\) on \(Q\).  The induced map on
every derived-functor cohomology group is therefore zero.
\end{proof}

\subsection{Degree two: gerbes and the Brauer boundary}

\begin{definition}
\label{def:cohomological-brauer}
The fpqc cohomological Brauer group of a semiring scheme \(X\) is
\[
 \Br'_{\fpqc}(X)=H^2_{\fpqc}(X,\Gm).
\]
This notation makes no assertion that the group is represented by
Azumaya semialgebras.
\end{definition}

The notation is used only for the right-hand edge of the
universal-coefficient filtration.  None of the results below identifies
\(\Br'_{\fpqc}(X)\) with an Azumaya Brauer group or computes it.  In
particular, the root-gerbe subgroup constructed below exists and its
order is determined before any such comparison is available.

\begin{theorem}[Diagonalizable gerbe sequence]
\label{thm:diagonalizable-gerbe-sequence}
Equivalence classes of fpqc gerbes banded by \(D_X(Q)\) form
\(H^2_{\fpqc}(X,D_X(Q))\), and there is a natural exact sequence
\begin{equation}
\label{eq:diagonalizable-gerbe-sequence}
\begin{aligned}
 0&\longrightarrow \Ext^1_{\mathbb Z}(Q,\Pic(X))
 \longrightarrow H^2_{\fpqc}(X,D_X(Q))\\
 &\longrightarrow \Hom_{\mathbb Z}(Q,\Br'_{\fpqc}(X))
 \longrightarrow0.
\end{aligned}
\end{equation}
If \(Q=Q_M=\coker M\), this becomes
\begin{equation}
\label{eq:matrix-gerbe-sequence}
 0\longrightarrow
 \coker\Phi_M(\Pic(X))
 \longrightarrow
 H^2_{\fpqc}(X,K_M)
 \longrightarrow
 \ker\Phi_M(\Br'_{\fpqc}(X))
 \longrightarrow0.
\end{equation}
\end{theorem}

\begin{proof}
The classification of abelian banded gerbes by second cohomology is
the theorem of Giraud; see also
\cite{StacksSecondCohomologyGerbes}.  Apply
Theorem~\ref{thm:universal-coefficient} with \(i=2\), and then use
Lemma~\ref{lem:character-resolution} for the matrix form.
\end{proof}

For a line bundle \(L\), let
\(\sqrt[m]{L/X}\) be the stack whose objects over \(T\to X\) are
pairs
\[
 (N,\psi),\qquad
 N\ \text{a line bundle on }T,\qquad
 \psi\colon N^{\ot m}\longiso L_T.
\]
It is locally nonempty, and two objects differ by a
\(\mu_m\)-torsor; hence it is a \(\mu_m\)-banded fpqc gerbe.

\begin{theorem}[Root-gerbe detection]
\label{thm:root-gerbe-detection}
For every \(m\geq1\), the cyclic specialization of
\eqref{eq:diagonalizable-gerbe-sequence} is
\begin{equation}
\label{eq:cyclic-gerbe-sequence}
 0\longrightarrow
 \Pic(X)/m\Pic(X)
 \xrightarrow{\ \partial_m\ }
 H^2_{\fpqc}(X,\mu_m)
 \longrightarrow
 \Br'_{\fpqc}(X)[m]
 \longrightarrow0.
\end{equation}
The class \(\partial_m([L])\) is represented by
\(\sqrt[m]{L/X}\).  Consequently:
\begin{enumerate}
\item \(\sqrt[m]{L/X}\) is neutral if and only if
\([L]\in m\Pic(X)\);
\item the order of its class is the order of \([L]\) in
\(\Pic(X)/m\Pic(X)\).
\end{enumerate}
\end{theorem}

\begin{proof}
For \(Q=C_m\),
\[
\begin{aligned}
 \Ext^1_{\mathbb Z}(C_m,\Pic(X))
 &\cong\Pic(X)/m\Pic(X),\\
 \Hom_{\mathbb Z}(C_m,\Br'_{\fpqc}(X))
 &=\Br'_{\fpqc}(X)[m].
\end{aligned}
\]
The connecting map
\(\Pic(X)\to H^2_{\fpqc}(X,\mu_m)\) in the long exact sequence of
\[
 1\longrightarrow\mu_m\longrightarrow\Gm
 \xrightarrow{[m]}\Gm\longrightarrow1
\]
is the gerbe of lifts of a \(\Gm\)-torsor through \([m]\).  Under
\(H^1(X,\Gm)=\Pic(X)\), that lift gerbe is precisely
\(\sqrt[m]{L/X}\).  Exactness proves neutrality, while injectivity of
the first arrow of \eqref{eq:cyclic-gerbe-sequence} proves the order
statement.
\end{proof}

\begin{remark}
\label{rem:vertical-not-derived-semimodules}
Theorem~\ref{thm:derived-character-formula} is a derived statement
about abelian fpqc sheaves.  It neither presupposes nor constructs a
derived category of semimodules.  The semiring-specific input is
geometric: Proposition~\ref{prop:matrix-finite-torsor} makes
\eqref{eq:character-sheaf-resolution} exact on the fpqc site, and
Theorem~\ref{thm:global-graded-descent} represents its degree-one
classes by finite locally free torsors.
\end{remark}

\section{Scalar extension and classical shadows}
\label{sec:base-change}

The character resolution, its torsors, and its derived cohomology are
compatible with pullback.  This provides a comparison with ordinary
scheme theory when a semiring has a nondegenerate ring completion,
and it also isolates classes which are intrinsically
characteristic-one.

\begin{theorem}[Base-change compatibility]
\label{thm:base-change}
Let \(f\colon X'\to X\) be a morphism of semiring schemes and let
\(Q\) be a finite abelian group.
\begin{enumerate}
\item There are canonical identifications
\[
 D_X(Q)\times_XX'\cong D_{X'}(Q),
\qquad
 K_{M,X}\times_XX'\cong K_{M,X'}.
\]
\item A Picard-strong grading
\(\mathcal B=\bigoplus_{q\in Q}\mathcal L_q\)
pulls back to
\[
 f^*\mathcal B
 \cong
 \bigoplus_{q\in Q}f^*\mathcal L_q
\]
with its pulled-back strong multiplication maps.
\item Pullback commutes with the equivalences
\[
 \Kum_M(X)
 \longiso\Tors_{\fpqc}(X,K_M)
 \longiso\operatorname{PicGr}_{Q_M}(X)
\]
and with the degree-one and degree-two connecting morphisms.
\item The derived character formula is natural in the commutative
square
\[
\begin{CD}
 R\Gamma_{\fpqc}(X,D_X(Q))
 @>{\sim}>>
 R\Hom_{\mathbb Z}(Q,R\Gamma_{\fpqc}(X,\Gm))\\
 @VV{\operatorname{res}_f}V
 @VV{R\Hom(Q,\operatorname{res}_f)}V\\
 R\Gamma_{\fpqc}(X',D_{X'}(Q))
 @>{\sim}>>
 R\Hom_{\mathbb Z}(Q,R\Gamma_{\fpqc}(X',\Gm)).
\end{CD}
\]
\end{enumerate}
\end{theorem}

\begin{proof}
The monoid-semiring identity
\[
 A'\ot_AA[Q]\cong A'[Q]
\]
proves the first assertion on affine opens.  Pullback preserves
finite direct sums, tensor products, duals, and isomorphisms, proving
the second.  These operations also define the three functors in the
matrix Kummer equivalence, so the third assertion follows.  Finally,
pullback takes the character-sheaf resolution
\eqref{eq:character-sheaf-resolution} on \(X\) to the corresponding
resolution on \(X'\).  The vertical arrows in the displayed square
are the restriction maps induced by the adjunction morphism to the
derived direct image along
\(X'_{\fpqc}\to X_{\fpqc}\).  Applying this construction to both
two-term resolutions gives the same morphism of homotopy fibers;
this proves the last assertion and specifies its naturality.
\end{proof}

\subsection{Ring completion}

Let
\[
 A^{\mathrm{gp}}=\mathbb Z\ot_{\mathbb N}A
\]
be the universal ring receiving a semiring morphism from \(A\).  For
an \(A\)-semiring scheme \(X\), write
\[
 X_{\mathbb Z}
 =
 X\times_{\Spec A}\Spec A^{\mathrm{gp}}.
\]
The ring \(A^{\mathrm{gp}}\) may be zero.

\begin{proposition}[Classical shadow]
\label{prop:ring-completion}
Every finite diagonalizable torsor, Picard-strong grading, matrix
Kummer object, and diagonalizable gerbe on \(X\) has a canonical
pullback, or classical shadow, on \(X_{\mathbb Z}\).  The universal torus cover
\[
 \varphi_M\colon(\Gm{}_{\mathbb N})^r
 \longrightarrow(\Gm{}_{\mathbb N})^r
\]
becomes the classical finite locally free torus isogeny over
\(\mathbb Z\).  The restriction maps for diagonalizable and
\(\Gm\)-cohomology form the natural commutative square of
Theorem~\ref{thm:base-change}; no cohomological base-change
isomorphism is asserted.
\end{proposition}

\begin{proof}
Apply Theorem~\ref{thm:base-change} to the morphism
\(X_{\mathbb Z}\to X\) induced by \(A\to A^{\mathrm{gp}}\).
This gives the asserted pullbacks of torsors, gradings, Kummer objects,
and gerbes, as well as the natural restriction square for the derived
character formula.  For the universal cover, specialize to
\(A=\mathbb N\) and \(A^{\mathrm{gp}}=\mathbb Z\).  Then the Hopf
semiring \(\mathbb N[Q]\) becomes the usual Hopf algebra
\(\mathbb Z[Q]\), and \(\mathbb N[\mathbb Z^r]\) becomes
\(\mathbb Z[\mathbb Z^r]\).  Proposition~\ref{prop:matrix-finite-torsor}
therefore base-changes to the classical finite locally free torus
isogeny over \(\mathbb Z\).
\end{proof}

\begin{remark}
\label{rem:nonconservative-shadow}
The classical-shadow functor is not conservative.  If \(A\) is a
nonzero idempotent semiring, then \(A^{\mathrm{gp}}=0\): the relation
\(1+1=1\) becomes \(1=0\) after additive group completion.  The
Boolean root gerbes of
Corollary~\ref{cor:boolean-cohomological-jump} therefore record
cohomology which cannot be reconstructed from an ordinary ring
fiber.
\end{remark}

\section{Torsor--gerbe separation in characteristic one}
\label{sec:tropical}

The universal-coefficient filtration detects a phenomenon which
degree-one Kummer theory cannot see: a semiring scheme may have no
nontrivial finite diagonalizable torsors of a given type while
carrying nontrivial gerbes of that type.  Tropical projective space
provides a uniform family.
These calculations are deliberately consequences of Jun's Picard
group theorem, not new computations of the Picard or cohomological
Brauer groups.  Their purpose is to show that the two edges of the
character filtration have geometrically different behavior.  The
independent affine example in
Section~\ref{sec:affine-boolean-torsor} supplies a nonzero degree-one
class with an explicit coordinate semialgebra.

\subsection{Projective cohomology in degrees one and two}

Let \(S\) be a totally ordered idempotent semifield and \(d\geq1\).
For the projective semiring scheme \(\mathbb P^d_S\), Jun proves
\begin{equation}
\label{eq:jun-projective-computation}
 \Gamma(\mathbb P^d_S,\OO)=S,
 \qquad
 \Pic(\mathbb P^d_S)\cong\mathbb Z,
\end{equation}
with \(\OO(1)\) as a generator
\cite[Propositions~4.9 and~4.25 and
Corollary~4.29]{Jun2017}.  We use the semiring-scheme model fixed in
Section~\ref{sec:foundations}; related tropical frameworks are
developed in
\cite{GiansiracusaGiansiracusa2016,MaclaganRincon2018,
Lorscheid2022,Lorscheid2023}.

\begin{theorem}[Projective diagonalizable filtration]
\label{thm:projective-diagonalizable-filtration}
Let \(Q\) be a finite abelian group.  There are natural exact
sequences
\begin{equation}
\label{eq:projective-H1}
 0\longrightarrow
 \Ext^1_{\mathbb Z}(Q,S^\times)
 \longrightarrow
 H^1_{\fpqc}(\mathbb P^d_S,D(Q))
 \longrightarrow0
\end{equation}
and
\begin{equation}
\label{eq:projective-H2}
 0\longrightarrow
 \Ext^1_{\mathbb Z}(Q,\mathbb Z)
 \longrightarrow
 H^2_{\fpqc}(\mathbb P^d_S,D(Q))
 \longrightarrow
 \Hom_{\mathbb Z}\!\left(
 Q,\Br'_{\fpqc}(\mathbb P^d_S)
 \right)
 \longrightarrow0.
\end{equation}
In particular, the left subgroup in degree two is finite of order
\(|Q|\).
\end{theorem}

\begin{proof}
Because \(Q\) is finite and \(\mathbb Z\) is torsion free,
\(\Hom_{\mathbb Z}(Q,\mathbb Z)=0\).  Substitute
\eqref{eq:jun-projective-computation} in
\eqref{eq:invariant-degree-one} and
\eqref{eq:diagonalizable-gerbe-sequence}.  Finally, if
\[
 Q\cong\bigoplus_j\mathbb Z/d_j\mathbb Z,
\]
then
\[
 \Ext^1_{\mathbb Z}(Q,\mathbb Z)
 \cong
 \bigoplus_j\mathbb Z/d_j\mathbb Z
\]
after choosing invariant-factor generators, so its order is
\(\prod_jd_j=|Q|\).  The group is intrinsically the Pontryagin dual
\(\Hom(Q,\mathbb Q/\mathbb Z)\), although the displayed
identification with \(Q\) is not canonical.
\end{proof}

\begin{corollary}[Matrix form]
\label{cor:projective-matrix-filtration}
Let \(M\in\operatorname{Mat}_{r\times r}(\mathbb Z)\) have nonzero
determinant and \(K_M=D(Q_M)\).  Then
\[
 H^1_{\fpqc}(\mathbb P^d_S,K_M)
 \cong
 (S^\times)^r/\Phi_M((S^\times)^r),
\]
whereas degree two sits in
\begin{equation}
\label{eq:projective-matrix-H2}
\begin{aligned}
 0&\longrightarrow
 \coker\!\left(
 M^{\mathsf T}\colon\mathbb Z^r\to\mathbb Z^r\right)
 \longrightarrow H^2_{\fpqc}(\mathbb P^d_S,K_M)\\
 &\longrightarrow
 \ker\Phi_M\!\left(\Br'_{\fpqc}(\mathbb P^d_S)\right)
 \longrightarrow0.
\end{aligned}
\end{equation}
The canonical left subgroup has order \(|\det M|\).
\end{corollary}

\begin{proof}
Use Lemma~\ref{lem:character-resolution} and note that
\(\Phi_M\) acts as \(M^{\mathsf T}\) on
\(\Pic(\mathbb P^d_S)^r\cong\mathbb Z^r\).
\end{proof}

\subsection{Vanishing torsors and surviving root gerbes}

Use the max-plus semifields
\[
 \TT_{\mathbb R}=(\mathbb R\cup\{-\infty\},\max,+),
 \qquad
 \TT_{\mathbb Z}=(\mathbb Z\cup\{-\infty\},\max,+),
\]
and the Boolean semifield \(\mathbb B\).  Their unit groups are,
respectively, \((\mathbb R,+)\), \((\mathbb Z,+)\), and the trivial
group.

\begin{theorem}[Torsor--gerbe separation]
\label{thm:torsor-gerbe-separation}
Let \(m\geq2\) and \(d\geq1\).  For
\(S=\TT_{\mathbb R}\) or \(S=\mathbb B\),
\begin{equation}
\label{eq:torsor-vanishing-gerbe-survival}
 H^1_{\fpqc}(\mathbb P^d_S,\mu_m)=0,
 \qquad
 \mathbb Z/m\mathbb Z
 \lhook\joinrel\longrightarrow
 H^2_{\fpqc}(\mathbb P^d_S,\mu_m).
\end{equation}
The generator of the displayed subgroup is the root gerbe
\[
 \sqrt[m]{\OO(1)/\mathbb P^d_S},
\]
and it has exact order \(m\).  For the discrete tropical semifield,
\[
 H^1_{\fpqc}(\mathbb P^d_{\TT_{\mathbb Z}},\mu_m)
 \cong\mathbb Z/m\mathbb Z,
\]
while the same degree-two root-gerbe subgroup survives.
\end{theorem}

\begin{proof}
For \(\TT_{\mathbb R}\), multiplication by \(m\) on the unit group
\(\mathbb R\) is bijective; for \(\mathbb B\), the unit group is
trivial.  Since
\(\Pic(\mathbb P^d_S)\cong\mathbb Z\) has no \(m\)-torsion, the
degree-one cyclic Kummer sequence gives the first vanishing in
\eqref{eq:torsor-vanishing-gerbe-survival}.  Over
\(\TT_{\mathbb Z}\), its unit quotient is
\(\mathbb Z/m\mathbb Z\), giving the displayed degree-one group.

In every case,
\[
 \Pic(\mathbb P^d_S)/m\Pic(\mathbb P^d_S)
 \cong\mathbb Z/m\mathbb Z,
\]
generated by \([\OO(1)]\).  Theorem
\ref{thm:root-gerbe-detection} injects this quotient into
\(H^2_{\fpqc}(\mathbb P^d_S,\mu_m)\) and identifies its generator
with the stated root gerbe.  Its order is therefore exactly \(m\).
\end{proof}

\begin{corollary}[Boolean cohomological jump]
\label{cor:boolean-cohomological-jump}
On \(\mathbb P^d_{\mathbb B}\), every \(\mu_m\)-torsor is trivial,
but the \(\mu_m\)-banded gerbe of \(m\)th roots of \(\OO(1)\) is
nonneutral.  This class has no nondegenerate classical shadow under
ring completion.
\end{corollary}

\begin{proof}
Only the last assertion remains.  Additive group completion of a
nonzero idempotent semiring is the zero ring, because \(1+1=1\)
becomes \(1=0\).  Thus the classical base change of
\(\mathbb P^d_{\mathbb B}\) is degenerate, while the semiring
root-gerbe class is nonzero by
Theorem~\ref{thm:torsor-gerbe-separation}.
\end{proof}

\begin{remark}
\label{rem:degree-separation}
The distinction is structural, not a failure to find enough
base-valued roots.  Degree one measures power classes of global units
and torsion in \(\Pic\); degree two receives the quotient
\(\Pic/m\Pic\).  For a torsion-free but non-\(m\)-divisible Picard
group, the first Picard contribution vanishes while the second is
forced to survive.
\end{remark}

\section{Concluding remarks}
\label{sec:conclusion}

The paper has two structural layers.  The first is geometric.
A finite diagonalizable torsor over a semiring scheme is not treated
as an undifferentiated fpqc form of a free semimodule of rank
\(|Q|\).  Its coaction splits the coordinate algebra into
homogeneous fpqc forms of rank one, and the torsor identity forces
the multiplication maps among those pieces to be isomorphisms.
Rank-one descent turns them into line bundles.  Their finite direct
sum then proves Zariski-local freeness of the entire torsor.  This
proves local freeness for diagonalizable torsors without asserting
the still-unknown higher-rank implication for general semimodules.

The second layer records cohomological consequences.  Once the finite
torus resolution of \(D_X(Q)\) is known to be fpqc-exact, standard
derived internal Hom gives, for every \(f\colon X\to Y\),
\[
 Rf_*D_X(Q)
 \simeq
 R\mathcal{H}\!om_{\mathbb Z,Y}
 \bigl(Q_Y,Rf_*\Gm{}_{X}\bigr);
\]
after derived global sections it becomes
\[
 R\Gamma_{\fpqc}(X,D_X(Q))
 \simeq
 R\Hom_{\mathbb Z}
 \bigl(Q,R\Gamma_{\fpqc}(X,\Gm)\bigr).
\]
The resulting universal-coefficient filtration is intrinsic to the
character group and independent of a matrix presentation.  In
degree one it explains the two components of Kummer torsors:
\(\Ext^1(Q,\Gamma(\OO_X^\times))\) and
\(\Hom(Q,\Pic(X))\).  In degree two the same mechanism replaces
units by the Picard group and Picard classes by the cohomological
Brauer group.  Thus matrix Kummer torsors, diagonalizable gerbes, and
root gerbes appear as successive edges of one character resolution.
This organizing statement reduces their cohomology to
\(\Gm\)-cohomology; it does not independently compute the latter.

The degree-one edge has two concrete geometric consequences.  Every
homomorphism \(Q\to\Pic(X)\) is realized by a finite locally free
diagonalizable torsor, and the ambiguity in its coherent
multiplication is measured by
\(\Ext^1_{\mathbb Z}(Q,\Gamma(X,\OO_X^\times))\).  On the affine
Boolean torus \(X=\operatorname{Spec}\mathbb B[u^{\pm1}]\), this unit
term supplies an explicit subgroup \(\mathbb Z/m\mathbb Z\) of
\(H^1_{\fpqc}(X,\mu_m)\), represented by the finite free root
semialgebras \(z^m=u^a\).  Thus the geometric theory has nontrivial
affine classes even when no projective Picard computation is used.

The projective idempotent examples show why the higher layer is
necessary.  On \(\mathbb P^d_{\mathbb B}\) and
\(\mathbb P^d_{\mathbb T_{\mathbb R}}\), degree-one
\(\mu_m\)-cohomology vanishes, but the degree-two boundary of
\(\OO(1)\) has exact order \(m\).  Over the Boolean semifield its
ring completion is zero.  These gerbes are therefore classes not
detected by passage to the ordinary ring-completion fiber.

\subsection{Four concrete problems}

The limitations above suggest the following questions.  We state them
in deliberately small forms: even these first cases appear to require
ideas beyond the character-by-character argument used here.

\begin{enumerate}
\item \emph{Nonsplit multiplicative type.}
Let \(G\) be a form of \(\mu_m\) that becomes diagonalizable after a
finite faithfully flat extension.  Does every fpqc \(G\)-torsor over a
semiring scheme become finite locally free in the Zariski topology?
The first test is a nontrivial quadratic form of \(\mu_m\).  After the
splitting extension the character pieces are line bundles, but their
descent data can permute the pieces, so the present proof no longer
descends them one at a time.

\item \emph{Equivariant higher-rank descent.}
Which finite flat group schemes have the property that every torsor's
regular coaction decomposes its coordinate semimodule into summands to
which rank-one descent applies?  A useful intermediate goal is to
classify finite group actions whose isotypic pieces are invertible
semimodules.  Such a criterion would locate the exact boundary between
the diagonalizable theorem proved here and the unresolved general
higher-rank problem.

\item \emph{Azumaya classes versus cohomological gerbes.}
Is there a natural Azumaya Brauer group of semiring schemes together
with a comparison map to \(H^2_{\fpqc}(X,\Gm)\)?  Is that map injective
or surjective for either
\(X=\operatorname{Spec}\mathbb B[u^{\pm1}]\) or
\(X=\mathbb P^1_{\mathbb B}\)?  The Boolean root gerbes constructed
here provide explicit classes against which any proposed comparison
can be tested.

\item \emph{Changing the tropical model.}
Does the \(\mu_m\)-root gerbe of \(\OO(1)\) on
\(\mathbb P^1_{\mathbb B}\) persist under comparison with congruence,
hyperfield, or valuated-matroid models of tropical geometry?  The
smallest case \(d=1\) is already meaningful: the class is nonzero in
the relative semiring-scheme model although its ordinary ring
completion is degenerate.
\end{enumerate}

These problems also explain the role of the derived formula.  It does
not finish the cohomology theory; rather, it identifies precisely which
unit, Picard, and Brauer-type inputs a future theory must compute or
compare.

\section*{Acknowledgments and disclosures}
The authors received no external funding for this work and declare no
competing interests.  No research data were used.

\bibliographystyle{amsalpha}
\bibliography{references}

\end{document}